\documentclass{article}
\usepackage{amssymb}

\input{tcilatex}
\begin{document}

\begin{center}
\textbf{Nonlocal Cauchy problems for wave equations and applications}

\textbf{Veli\ B. Shakhmurov}

Department of Mechanical Engineering, Okan University, Akfirat, Tuzla 34959
Istanbul, Turkey,

E-mail: veli.sahmurov@okan.edu.tr;

\textbf{Abstract}
\end{center}

In this paper, the existence, the uniqueness and estimates of solution to
the integral Cauchy problem for linear and nonlinear abstract wave equations
are proved. The equation includes a linear operator $A$ defined in a Banach
space $E$, in which by choosing $E$ and $A$ we can obtain numerous classis
of nonlocal initial value problems for wave equations which occur in a wide
variety of physical systems.

\textbf{Key Word:}$\mathbb{\ \ }$wave equations\textbf{, }Semigroups of
operators, Hyperbolic-operator equations, cosine operator functions,
operator valued $L^{p}-$Fourier multipliers

\begin{center}
\bigskip\ \ \textbf{AMS: 65J, 65N, 35J, 47D}

\textbf{1}. \textbf{Introduction}
\end{center}

The subject of this paper is to study the existence, uniqueness and
regularity properties to solution of the integral Cauchy problem (ICP) for
the nonlinear abstract wave equation (NAWE)%
\begin{equation}
u_{tt}-\Delta u+Au+F\left( u\right) =0,\text{ }x\in R^{n},\text{ }t\in
\left( 0,\infty \right) ,  \tag{1.1}
\end{equation}%
\begin{equation}
u\left( x,0\right) =\varphi \left( x\right) +\dint\limits_{0}^{T}\alpha
\left( \sigma \right) u\left( x,\sigma \right) d\sigma ,\text{ }  \tag{1.2}
\end{equation}

\[
u_{t}\left( x,0\right) =\psi \left( x\right) +\dint\limits_{0}^{T}\beta
\left( \sigma \right) u_{t}\left( x,\sigma \right) d\sigma , 
\]%
where $A$ is a linear operator in a Banach space $E,$ $\alpha \left(
s\right) $ and $\beta \left( s\right) $ are measurable functions on $\left(
0,T\right) $, $u(x,t)$ denotes the $E$-valued unknown function, $f(u)$ is
the given nonlinear function, $\varphi \left( x\right) $ and $\psi \left(
x\right) $ are the given initial value functions, subscript $t$ indicates
the partial derivative with respect to $t$, $n$ is the dimension of space
variable $x$ and $\Delta $ denotes the Laplace operator in $R^{n}.$

\textbf{Remark 1.1. }Note that, particularly, the conditions $\left(
1.2\right) $\ can be expressed as the following multipoint nonlocal
conditions 
\begin{equation}
u\left( x,0\right) =\varphi \left( x\right) +\dsum\limits_{k=1}^{l}\alpha
_{k}u\left( x,\lambda _{k}\right) ,\text{ }u_{t}\left( 0,x\right) =\psi
\left( x\right) +\dsum\limits_{k=1}^{m}\beta _{k}u_{t}\left( x,\lambda
_{k}\right) ,  \tag{1.3}
\end{equation}%
where $l$ is a positive integer, $\alpha _{k}$, $\beta _{k}$ are complex
numbers and $\lambda _{k}\in \left( 0,\infty \right) .$

Since the Banach space $E$ is arbitrary and $A$ is a possible linear
operator, by choosing $E$, $A$ and integral conditions, we can obtain
numerous classis of nonlocal initial value problems for wave equations which
occur in a wide variety of physical systems, particularly in the propagation
of longitudinal deformation waves in an elastic rod, hydro-dynamical process
in plasma, in materials science which describe spinodal decomposition, in
the absence of mechanical stresses (see $\left[ 29\right] $, $\left[ 45%
\right] $, $\left[ 48\right] $ and the references therein). If $F\left(
u\right) =$ $\lambda \left\vert u\right\vert ^{p}u,$ from $\left( 1.1\right) 
$ we get the ICP for the following NAWE%
\begin{equation}
u_{tt}-\Delta u+Au+\lambda \left\vert u\right\vert ^{p}u=0,\text{ }x\in
R^{n},\text{ }t\in \left( 0,\infty \right) ,  \tag{1.4}
\end{equation}%
\[
u\left( x,0\right) =\varphi \left( x\right) +\dint\limits_{0}^{T}\alpha
\left( \sigma \right) u\left( x,\sigma \right) d\sigma ,\text{ } 
\]

\[
u_{t}\left( x,0\right) =\psi \left( x\right) +\dint\limits_{0}^{T}\beta
\left( \sigma \right) u_{t}\left( x,\sigma \right) d\sigma , 
\]%
where $p\in \left( 1,\infty \right) $, $\lambda $ is a real number,

Let $\mathbb{N},$ $\mathbb{R}$ and $\mathbb{C}$\ denote the sets of all
natural, real and complex numbers, respectively. For $E=\mathbb{C}$, $\alpha
_{k}=\beta _{k}=0$ and $A=0$ the problem $\left( 1.4\right) $\ become as the
following integral Cauchy problem 
\begin{equation}
u_{tt}-\Delta u+\lambda \left\vert u\right\vert ^{p-1}u=0,\text{ }x\in R^{n},%
\text{ }t\in \left( 0,\infty \right) ,  \tag{1.5}
\end{equation}%
\[
u\left( x,0\right) =\varphi \left( x\right) +\dint\limits_{0}^{T}\alpha
\left( \sigma \right) u\left( x,\sigma \right) d\sigma ,\text{ } 
\]

\[
u_{t}\left( x,0\right) =\psi \left( x\right) +\dint\limits_{0}^{T}\beta
\left( \sigma \right) u_{t}\left( x,\sigma \right) d\sigma , 
\]

If we\ choose $E$ a concrete space, for example $E=L^{2}\left( \Omega
\right) $, $A=L,$ where $\Omega $ is a domain in $R^{m}$ with sufficiently
smooth boundary and $L$ is an elliptic operator in $L^{2}\left( \Omega
\right) ,$ then from $\left( 1.1\right) -\left( 1.2\right) $ we obtain the
existence, uniqueness and the regularity properties of the mixed problem for
linear wave equation

\[
u_{tt}-\Delta u+Lu=F\left( x,t\right) ,\text{ }t\in \left[ 0,T\right] \text{%
, }x\in R^{n},\text{ }y\in \Omega , 
\]%
and the following nonlinear wave equations (NWE) equation%
\[
u_{tt}-\Delta u+Lu+\lambda \left\vert u\right\vert ^{p-1}u=0,\text{ }t\in %
\left[ 0,T\right] \text{, }x\in R^{n},\text{ }y\in \Omega , 
\]%
where $u=u\left( x,y,t\right) .$

\ Moreover, let we choose $E=L^{2}\left( 0,1\right) $ and $A$ to be
differential operator with generalized Wentzell-Robin boundary condition
defined by 
\[
D\left( A\right) =\left\{ u\in W^{2,2}\left( 0,1\right) ,\text{ }%
B_{j}u=Au\left( j\right) =0,\text{ }j=0,1\right\} ,\text{ } 
\]%
\begin{equation}
\text{ }Au=au^{\left( 2\right) }+bu^{\left( 1\right) }  \tag{1.6}
\end{equation}%
where $a=a\left( y\right) $ and $b=b\left( y\right) $ are complex-valued
functions. Then, from the main our theorem we get the existence, uniqueness
and regularity properties of multipoint Wentzell-Robin type mixed problem
for the linear wave equation 
\begin{equation}
u_{tt}-\Delta u+au_{yy}+bu_{y}=F\left( x,y,t\right) ,\text{ }  \tag{1.7}
\end{equation}%
\ \ \ 

\begin{equation}
u\left( x,y,0\right) =\varphi \left( x,y\right) +\dint\limits_{0}^{T}\alpha
\left( \sigma \right) u\left( x,y,\sigma \right) d\sigma ,\text{ }  \tag{1.8}
\end{equation}

\[
u_{t}\left( x,y,0\right) =\psi \left( x,y\right) +\dint\limits_{0}^{T}\beta
\left( \sigma \right) u_{t}\left( x,y,\sigma \right) d\sigma , 
\]

\begin{equation}
a\left( j\right) u_{yy}\left( x,j,t\right) +b\left( j\right) u_{y}\left(
x,j,t\right) =0,\text{ }j=0,1\text{ for all }t\in \left[ 0,T\right] 
\tag{1.9}
\end{equation}%
and for the following NWE%
\begin{equation}
u_{tt}-\Delta u+au_{yy}+bu_{y}+F\left( u\right) =0,\text{ }  \tag{1.10}
\end{equation}%
where 
\[
u=u\left( x,y,t\right) \text{, }x\in R^{n},\text{ }y\in \left( 0,1\right) ,%
\text{ }t\in \left[ 0,T\right] . 
\]%
Note that, the regularity properties of Wentzell-Robin type boundary value
problems (BVPs) for elliptic equations were studied e.g. in $\left[ \text{%
12, 25 }\right] $ and the references therein. Moreover, if put $E=l_{2}$ and 
$A$ choose as a infinite matrix $\left[ a_{mj}\right] $, $m,j=1,2,...,\infty
,$ then from our results we obtain the existence, uniqueness and regularity
properties of integral Cauchy problem for infinity many system of linear
wave equations 
\begin{equation}
\left( u_{m}\right) _{tt}-\Delta
u_{m}+\sum\limits_{j=1}^{N}a_{mj}u_{j}=F_{j}\left( x,t\right) ,\text{ }x\in
R^{n},\text{ }t\in \left[ 0,T\right] \text{,}  \tag{1.11}
\end{equation}%
\[
u_{m}\left( x,0\right) =\varphi _{m}\left( x\right)
+\dint\limits_{0}^{T}\alpha \left( \sigma \right) u_{m}\left( x,\sigma
\right) d\sigma ,\text{ } 
\]

\[
\left( u_{m}\right) _{t}\left( x,0\right) =\psi \left( x\right)
+\dint\limits_{0}^{T}\beta \left( \sigma \right) \left( u_{m}\right)
_{t}\left( x,\sigma \right) d\sigma , 
\]%
and infinity many system of NWE equation 
\begin{equation}
\left( u_{m}\right) _{tt}-\Delta
u_{m}+\sum\limits_{j=1}^{N}a_{mj}u_{j}+F_{m}\left(
u_{1},u_{2},...u_{N}\right) =0,\text{ }x\in R^{n},\text{ }t\in \left[ 0,T%
\right] \text{,}  \tag{1.12}
\end{equation}%
where $a_{mj}$ are complex numbers, $u_{j}=u_{j}\left( x,t\right) .$

The existence of solutions and regularity properties of Cauchy problem for
NWE studied e.g in $\left[ \text{5, 7}\right] $, $\left[ 14\right] $, $\left[
18\right] $, $\left[ 20-24\right] $, $\left[ 26\right] $, $\left[ 30\right] $%
, $\left[ 34\right] ,$ $\left[ \text{42-46}\right] $ and the references
therein.\ In contrast to the mentioned above results we will study the
regularity properties of the problem $\left( 1.1\right) -\left( 1.2\right) $%
. Abstract differential equations studied e.g. in $\left[ 1-3\right] $, $%
\left[ 7-11\right] $, $\left[ 13\right] $, $\left[ 16-17\right] $, $\left[ 19%
\right] $, $\left[ 28\right] $, $\left[ 32\text{, }33\right] $, $\left[ 35-39%
\right] ,$ $\left[ 44\right] $ and $\left[ 47\right] .$ The Cauchy problems
for abstract hyperbolic equations were treated e.g. in $\left[ 3\right] ,$ $%
\left[ 7\right] $, $\left[ 9\right] $, $\left[ 11\right] $, $\left[ 16\right]
$, $\left[ \text{32, 33}\right] $.

The strategy is to express the wave equation as an integral equation, to
treat in the nonlinearity as a small perturbation of the linear part of the
equation, then use the contraction mapping theorem and utilize an estimate
for solutions of the linearized version to obtain a priori estimates on $%
L^{p}$ norms of solutions. Harmonic analysis, the operator theory,
interpolation of Banach Spaces, embedding theorems in abstract Sobolev
spaces are the main tools implemented to carry out the analysis.

Sometimes we use one and the same symbol $C$ (or $M$) without distinction in
order to denote positive constants which may differ from each other even in
a single context. When we want to specify the dependence of such a constant
on a parameter, say $\alpha $, we write $C_{\alpha }$ (or $M_{\alpha }$).

In order to state our results precisely, we introduce some notations and
some function spaces:

\begin{center}
\textbf{Definitions and} \textbf{Background}
\end{center}

Let $E$ be a Banach space. $L^{p}\left( \Omega ;E\right) $ denotes the space
of strongly measurable $E$-valued functions that are defined on the
measurable subset $\Omega \subset R^{n}$ with the norm

\[
\left\Vert f\right\Vert _{L^{p}}=\left\Vert f\right\Vert _{L^{p}\left(
\Omega ;E\right) }=\left( \int\limits_{\Omega }\left\Vert f\left( x\right)
\right\Vert _{E}^{p}dx\right) ^{\frac{1}{p}},1\leq p<\infty ,\text{ } 
\]

\[
\left\Vert f\right\Vert _{L^{\infty }}\ =\text{ess}\sup\limits_{x\in \Omega
}\left\Vert f\left( x\right) \right\Vert _{E}. 
\]%
The Banach space\ $E$ is called an UMD-space if\ the Hilbert operator%
\[
\left( Hf\right) \left( x\right) =\lim\limits_{\varepsilon \rightarrow
0}\int\limits_{\left\vert x-y\right\vert >\varepsilon }\frac{f\left(
y\right) }{x-y}dy 
\]%
\ is bounded in $L^{p}\left( R,E\right) ,$ $p\in \left( 1,\infty \right) $
(see. e.g. $\left[ 4\right] $). Any Hilbert space is a UMD space. Moreover,
UMD spaces include e.g. $L^{p}$, $l_{p}$ spaces and Lorentz spaces $L_{pq},$ 
$p$, $q\in \left( 1,\infty \right) $.

Here, \ 
\[
S_{\psi }=\left\{ \lambda \in \mathbb{C}\text{, }\left\vert \arg \lambda
\right\vert \leq \omega ,\text{ }0\leq \omega <\pi \right\} , 
\]%
\[
S_{\omega ,\varkappa }=\left\{ \lambda \in S_{\omega }\text{, }\left\vert
\lambda \right\vert >\varkappa >0\right\} \text{ }. 
\]

A closed linear operator\ $A$ is said to be sectorial in a Banach\ space $E$
with bound $M>0$ if $D\left( A\right) $ and $R\left( A\right) $ are dense on 
$E,$ $N\left( A\right) =\left\{ 0\right\} $ and $\left\Vert \left( A+\lambda
I\right) ^{-1}\right\Vert _{B\left( E\right) }\leq M\left\vert \lambda
\right\vert ^{-1}$ for any $\lambda \in S_{\omega },$ $0\leq \omega <\pi ,$
where $I$ is the identity operator in $E,$ $B\left( E\right) $ is the space
of bounded linear operators in $E;$ $D\left( A\right) $ and $R\left(
A\right) $ denote domain and range of the operator $A.$ It is known that
(see e.g.$\left[ \text{40, \S 1.15.1}\right] $) there exist fractional
powers\ $A^{\theta }$ of a sectorial operator $A.$ Let $E\left( A^{\theta
}\right) $ denote the space $D\left( A^{\theta }\right) $ with the graphical
norm 
\[
\left\Vert u\right\Vert _{E\left( A^{\theta }\right) }=\left( \left\Vert
u\right\Vert ^{p}+\left\Vert A^{\theta }u\right\Vert ^{p}\right) ^{\frac{1}{p%
}},1\leq p<\infty ,\text{ }0<\theta <\infty . 
\]

A closed linear operator\ $A$ belong to $\sigma \left( M_{0},\omega
,E\right) $ (see $\left[ 11\right] $, \S\ 11.2) if $D\left( A\right) $ is
dense on $E,$ the resolvent $\left( A-\lambda ^{2}I\right) ^{-1}$ exists for 
$\func{Re}\lambda >\omega $ and 
\[
\left\Vert \left( A-\lambda ^{2}I\right) ^{-1}\right\Vert _{B\left( E\right)
}\leq M_{0}\left\vert \func{Re}\lambda -\omega \right\vert ^{-1}\text{. } 
\]

\textbf{Remark 1.2. }Let $0\leq \gamma <1$\textbf{\ }It is known\ that if $%
A\in \sigma \left( M_{0},\omega ,E\right) $, then it is is a sectorial
operator in $E$ and it is an infinitesimal generator of $C_{0}$\ group of
bounded linear operator $U_{A}\left( t\right) $ satisfying 
\begin{equation}
\left\Vert U_{A}\left( t\right) \right\Vert _{B\left( E\right) }\leq
M_{0}e^{\omega \left\vert t\right\vert },\text{ }t\in \left( -\infty ,\infty
\right) ,  \tag{1.1}
\end{equation}%
\[
\left\Vert A^{\gamma }U_{A}\left( t\right) \right\Vert _{B\left( E\right)
}\leq M_{0}\left\vert t\right\vert ^{-\gamma },\text{ }t\in \left( -\infty
,\infty \right) 
\]%
(see e.g. $\left[ \text{33}\right] ,$ $\left[ \text{\S\ 1.6}\right] $,
Theorem 6.3).

Let $E_{1}$ and $E_{2}$ be two Banach spaces. $\left( E_{1},E_{2}\right)
_{\theta ,p}$ denotes the interpolation spaces obtained from $\left\{
E_{1},E_{2}\right\} $ by $K$-method, where, $\theta \in \left( 0,1\right) ,$ 
$p\in \left[ 1,\infty \right] $ \ $\left[ \text{40, \S 1.3.2}\right] $.

The sectorial operator $A$ is said to be $R-$sectorial in a Banach space $E$
if the set $\left\{ \xi \left( A+\xi \right) ^{-1}\text{: }\xi \in S_{\omega
}\right\} $, $0\leq \omega <\pi $ is $R$-bounded (see e.g. $\left[ 8\right] $%
) .

Let $E_{0}$ and $E$ be two Banach spaces and $E_{0}$ is continuously and
densely embedded into $E$. Let $\Omega $ be a domain in $R^{n}$ and $m$ is a
positive integer.\ $W^{m,p}\left( \Omega ;E_{0},E\right) $ denotes the space
of all functions $u\in L^{p}\left( \Omega ;E_{0}\right) $ that have the
generalized derivatives $\frac{\partial ^{m}u}{\partial x_{k}^{m}}\in
L^{p}\left( \Omega ;E\right) $ for $1\leq p\leq \infty $ with the norm 
\[
\ \left\Vert u\right\Vert _{W^{m,p}\left( \Omega ;E_{0},E\right)
}=\left\Vert u\right\Vert _{L^{p}\left( \Omega ;E_{0}\right)
}+\sum\limits_{k=1}^{n}\left\Vert \frac{\partial ^{m}u}{\partial x_{k}^{m}}%
\right\Vert _{L^{p}\left( \Omega ;E\right) }<\infty .
\]%
\ \ $\ \ $ For $E_{0}=E$ the space $W^{m,p}\left( \Omega ;E_{0},E\right) $
denotes by $W^{m,p}\left( \Omega ;E\right) .$ Here, $H^{s,p}\left(
R^{n};E\right) $, $-\infty <s<\infty $ denotes the $E-$valued Sobolev space
of order $s$ which is defined as: 
\[
H^{s,p}=H^{s,p}\left( R^{n};E\right) =\left( I-\Delta \right) ^{-\frac{s}{2}%
}L^{p}\left( R^{n};E\right) 
\]%
with the norm 
\[
\left\Vert u\right\Vert _{H^{s,p}}=\left\Vert \left( I-\Delta \right) ^{%
\frac{s}{2}}u\right\Vert _{L^{p}\left( R^{n};E\right) }<\infty .
\]%
It clear that $H^{0,p}\left( R^{n};E\right) =L^{p}\left( R^{n};E\right) .$
It is known that if $E$ is a UMD space, then $H^{m,p}\left( R^{n};E\right)
=W^{m,p}\left( R^{n};E\right) $ for positive integer $m$ (see e.g. $\left[ 
\text{41, \S\ 15}\right] $)$.$ $H^{s,p}\left( R^{n};E_{0},E\right) $ denote
the Sobolev-Lions type space, i.e., 
\[
H^{s,p}\left( R^{n};E_{0},E\right) =\left\{ u\in H^{s,p}\left(
R^{n};E\right) \cap L^{p}\left( R^{n};E_{0}\right) ,\right. \text{ }
\]%
\[
\left. \left\Vert u\right\Vert _{H^{s,p}\left( R^{n};E_{0},E\right)
}=\left\Vert u\right\Vert _{L^{p}\left( R^{n};E_{0}\right) }+\left\Vert
u\right\Vert _{H^{s,p}\left( R^{n};E\right) }<\infty \right\} .
\]%
$S\left( R^{n};E\right) $ denotes the Schwartz class, i.e., the space of $E$%
-valued rapidly decreasing smooth functions on $R^{n},$ equipped with its
usual topology generated by seminorms. Here, $S^{\prime }\left(
R^{n};E\right) $ denotes the space of all\ continuous linear operators $%
L:S\left( R^{n};E\right) \rightarrow E,$ equipped with the bounded
convergence topology. Recall that $S\left( R^{n};E\right) $ is norm dense in 
$L^{p}\left( R^{n};E\right) $ when $1\leq p<\infty $

Let $F$ denotes the Fourier transform. $\Psi \in L^{\infty }\left(
R^{n};B\left( E\right) \right) $ is called a multiplier from\ $L^{p}\left(
R^{n};E_{1}\right) $ to $L^{q}\left( R^{n};E\right) $ if there exists a
positive constant $C$ such that 
\[
\left\Vert F^{-1}\Psi \left( \xi \right) Fu\right\Vert _{L^{q}\left(
R^{n};E\right) }\leq C\left\Vert u\right\Vert _{L^{p}\left( R^{n};E\right) }%
\text{ for all }u\in S\left( R^{n};E_{1}\right) .
\]

We denote the set of all multipliers fom\ $L^{p}\left( R^{n};E\right) $ to $%
L^{q}\left( R^{n};E_{2}\right) $ by $M_{p}^{q}\left( E\right) $. Let $\Phi
_{h}=\left\{ \Psi _{h}\in M_{p}^{q}\left( E\right) ,\text{ }h\in \sigma
\right\} $ denote a collection of multipliers depending on the parameter $h.$

We say that $W_{h}$ is a uniform collection of multipliers if there exists a
positive constant $M$ independent of $h\in \sigma $ such that

\[
\left\Vert F^{-1}\Psi _{h}Fu\right\Vert _{L^{q}\left( R^{n};E\right) }\leq
M\left\Vert u\right\Vert _{L^{p}\left( R^{n};E\right) }\ \ \ \ \ \ 
\]%
for all $u\in S\left( R^{n};E\right) $ and $h\in \sigma .$

\textbf{Definition 1.1. } Assume $E$ is a Banach space and $r\in \left[ 1,2%
\right] .$ Suppose there exists a positive constant $C_{0}=C_{0}\left(
r,E\right) $ so that 
\[
\left\Vert Fu\right\Vert _{L^{r^{\prime }}\left( R^{n};E\right) }\leq
C_{0}\left\Vert Fu\right\Vert _{L^{r}\left( R^{n};E\right) } 
\]%
for $\frac{1}{r}+\frac{1}{r^{\prime }}=1$ and each $u\in S\left(
R^{n};E\right) .$ Then $E$ is called Fourier type $r.$

\textbf{Remark 1.3. }The simple estimate

\[
\left\Vert Ff\left( x\right) \right\Vert _{E}\leq \left\Vert f\right\Vert
_{L^{1}\left( R^{n};E\right) } 
\]%
shows that each Banach space $E$ has Fourier type $1.$ Bourgain $\left[ 4%
\right] $ has shown that each $B-$convex Banach space (thus, in particular,
each UMD space) has some non-trivial Fourier type $p>1.$

In order to define $E$-valued Besov spaces we consider the dyadic-like
subsets $\left\{ J_{k}\right\} _{k=0}^{\infty },$ $\left\{ I_{k}\right\}
_{k=0}^{\infty }$ of $R^{n}$ and partition of unity $\left\{ \varphi
_{k}\right\} _{k=0}^{\infty }$. Let $1\leq r,q\leq \infty $ and $s\in 
\mathbb{R}.$ The Besov space $B_{q,r}^{s}\left( R^{n};E\right) $ is the
space of all $f\in S^{\prime }\left( R^{n};E\right) $ with the norm 
\[
\left\Vert f\right\Vert _{B_{q,r}^{s}\left( R^{n};E\right) }=\left\Vert
\left\{ 2^{ks}\left( \check{\varphi}_{k}\ast f\right) \right\}
_{k=0}^{\infty }\right\Vert _{l_{r}\left( L^{q}\left( R^{n};E\right) \right)
}= 
\]%
\[
\left\{ 
\begin{array}{c}
\left[ \dsum\limits_{k=0}^{\infty }2^{ksr}\left\Vert \check{\varphi}_{k}\ast
f\right\Vert _{L^{q}\left( R^{n};E\right) }^{r}\right] ^{\frac{1}{r}}<\infty 
\text{, if }1\leq r<\infty \\ 
\sup\limits_{k\in \mathbb{N}_{0}}\left[ \dsum\limits_{k=0}^{\infty
}2^{ks}\left\Vert \check{\varphi}_{k}\ast f\right\Vert _{L^{q}\left(
R^{n};E\right) }\right] <\infty \text{, if }r=\infty%
\end{array}%
\right. . 
\]%
$B_{q,r}^{s}\left( R^{n};E\right) $-together with the above norm, is a
Banach space. it can be shown that different choices of $\left\{ \varphi
_{k}\right\} $ lead to equivalent norms on $B_{q,r}^{s}\left( R^{n};E\right)
.$Here, $B_{q,r}^{s}\left( R^{n};E_{0},E\right) $ denotes the space $%
L^{q}\left( R^{n};E_{0}\right) \cap B_{q,r}^{s}\left( R^{n};E\right) $ with
the norm 
\[
\left\Vert u\right\Vert _{B_{q,r}^{s}\left( R^{n};E_{0},E\right)
}=\left\Vert u\right\Vert _{L^{q}\left( R^{n};E_{0}\right) }+\left\Vert
u\right\Vert _{B_{q,r}^{s}\left( R^{n};E\right) }<\infty . 
\]

Let the operator $A$ be a generator of a strongly continuous cosine operator
function in a Banach space $E$ defined by formula%
\[
C\left( t\right) =\frac{1}{2}\left( e^{itA^{\frac{1}{2}}}+e^{-itA^{\frac{1}{2%
}}}\right) 
\]%
(see $\left[ \text{11, \S 11.2, 11.4}\right] $, or $\left[ \text{16}\right]
, $ $\left[ \text{33, 34}\right] $ ). Then, from the definition of sine
operator-function $S\left( t\right) $ we have%
\[
S\left( t\right) u=\dint\limits_{0}^{t}C\left( \sigma \right) ud\sigma 
\]%
and it follows that 
\[
S\left( t\right) u=\frac{1}{2i}A^{-\frac{1}{2}}\left( e^{itA^{\frac{1}{2}%
}}-e^{-itA^{\frac{1}{2}}}\right) . 
\]

\textbf{Lemma 1.1. }Let 
\[
\left\vert 1+\dint\limits_{0}^{T}\alpha \left( \sigma \right) \beta \left(
\sigma \right) d\sigma \right\vert >\dint\limits_{0}^{T}\left( \left\vert
\alpha \left( \sigma \right) \right\vert +\left\vert \beta \left( \sigma
\right) \right\vert \right) d\sigma . 
\]%
Then the operator $O$ defined by%
\[
O=\left[ 1+\dint\limits_{0}^{T}\dint\limits_{0}^{T}\alpha \left( \sigma
\right) \beta \left( \tau \right) d\sigma d\tau \right] I-\dint%
\limits_{0}^{T}\left( \alpha \left( s\right) +\beta \left( s\right) \right)
C\left( s\right) ds 
\]%
has an inverse $O^{-1}$ and the following estimate is satisfied 
\[
\left\Vert O^{-1}\right\Vert _{B\left( E\right) }\leq \left[ \left\vert
1+\dint\limits_{0}^{T}\alpha \left( s\right) \beta \left( s\right)
ds\right\vert -\dint\limits_{0}^{T}\left( \left\vert \alpha \left( s\right)
\right\vert +\left\vert \beta \left( s\right) \right\vert \right) ds\right]
^{-1}. 
\]

The embedding theorems in vector valued spaces play a key role in the theory
of DOEs. For estimating lower order derivatives we use following embedding
theorem that is obtained from $\left[ \text{37, Theorem 1}\right] $:

\textbf{Theorem A}$_{1}$. Suppose the following conditions are satisfied:

(1) $E$ is a Banach space and\ $A$ is a sectorial operator in $E;$

(2)\ $\alpha =\left( \alpha _{1},\alpha _{2},...,\alpha _{n}\right) $ is a $%
n $-tuples of nonnegative integer number\ and $s$ is a positive number such
that

$\varkappa =\frac{1}{s}\left[ \left\vert \alpha \right\vert +n\left( \frac{1%
}{p_{1}}-\frac{1}{p_{2}}\right) \right] \leq 1,$ $0\leq \mu \leq 1-\varkappa 
$, $1\leq p_{1}\leq p_{2}\leq \infty ;$ $0<h\leq h_{0},$ where $h_{0}$ is a
fixed positive number.

Then the embedding $D^{\alpha }B_{p_{1},r}^{s}\left( R^{n};E\left( A\right)
,E\right) \subset L_{p_{2}}\left( R^{n};E\left( A^{1-\varkappa -\mu }\right)
\right) $ is continuous and for $u\in H^{s,p}\left( R^{n};E\left( A\right)
,E\right) $ the following uniform estimate holds 
\[
\left\Vert D^{\alpha }u\right\Vert _{L^{p_{2}}\left( R^{n};E\left(
A^{1-\varkappa -\mu }\right) \right) }\leq h^{\mu }\left\Vert u\right\Vert
_{B_{p_{1},r}^{s}\left( R^{n};E\left( A\right) ,E\right) }+h^{-\left( 1-\mu
\right) }\left\Vert u\right\Vert _{L^{p_{1}}\left( R^{n};E\right) }. 
\]

By using $\left[ \text{15, Theorem 4.3}\right] $ we obtain:

\textbf{Proposition A}$_{1}.$ Assume the Banach spaces $E_{1},$ $E_{2}$ have
Fourier type $r\in \left[ 1,2\right] $ and%
\[
\Psi _{h}\in B_{r,1}^{n\left( \frac{1}{r}+\frac{1}{p_{1}}-\frac{1}{p_{2}}%
\right) }\left( R^{n};B\left( E_{1},E_{2}\right) \right) . 
\]%
Then $\Psi _{h}$ is a uniformly bounded collection of Fourier multiplier
from $L^{p_{1}}\left( R^{n};E\right) $ to $L^{p_{2}}\left( R^{n};E\right) $
for $p_{1}\leq p_{2}$ with $p_{1}$, $p_{2}\in \left[ 1,\infty \right] .$

\textbf{Proof. }First, in a similar way as in $\left[ \text{15, Theorem 4.3}%
\right] $ we show that $\Psi _{h}$ is a uniformly bounded collection of
Fourier multiplier from $L^{p}\left( R^{n};E\right) $ to $L^{p}\left(
R^{n};E\right) .$ Moreover, by Theorem A$_{1}$ we get that, for $s\geq
n\left( \frac{1}{p_{1}}-\frac{1}{p_{2}}\right) $ the embedding $%
B_{p_{1,1}}^{s}\left( R^{n};E\right) \subset L^{p_{2}}\left( R^{n};E\right) $
is continuous. From these two fact we obtain the conclusion.

The paper is organized as follows: In Section 1, some definitions and
background are given. In Section 2, we obtain the existence of unique
solution and a priory estimates for solution of the linearized problem $%
(1.1) $-$\left( 1.2\right) .$ In Section 3, we show the existence,
uniqueness and estimates of strong solution of the problem $(1.1)$-$\left(
1.2\right) $. In Section 4, the existence, uniqueness and a priory estimates
to solution of ICP for finite and infinite many system of wave equation is
derived.

\begin{center}
\textbf{2. Estimates for linearized equation}
\end{center}

In this section, we make the necessary estimates for solutions of ICP for
the abstract linear wave equation%
\begin{equation}
u_{tt}-\Delta u+Au=g\left( x,t\right) ,\text{ }x\in R^{n},\text{ }t\in
\left( 0,\infty \right) ,  \tag{2.1}
\end{equation}%
\begin{equation}
u\left( 0,x\right) =\varphi \left( x\right) +\dint\limits_{0}^{T}\alpha
\left( \sigma \right) u\left( \sigma ,x\right) d\sigma ,\text{ }  \tag{2.2}
\end{equation}

\[
u_{t}\left( 0,x\right) =\psi \left( x\right) +\dint\limits_{0}^{T}\beta
\left( \sigma \right) u_{t}\left( \sigma ,x\right) d\sigma . 
\]

\bigskip \textbf{Condition 2.1. }Assume:

(1) $\left\vert 1+\dint\limits_{0}^{T}\alpha \left( \sigma \right) \beta
\left( \sigma \right) d\sigma \right\vert >\dint\limits_{0}^{T}\left(
\left\vert \alpha \left( \sigma \right) \right\vert +\left\vert \beta \left(
\sigma \right) \right\vert \right) d\sigma ;$

(2) $E$ is a Banach space of Fourier type $r\in \left[ 1,2\right] $;

(3) $A\in \sigma \left( M_{0},\omega ,E\right) $ and $s>n\left( \frac{1}{r}+%
\frac{1}{p}\right) $ for $p\in \left[ 1,\infty \right] .$

Let 
\[
X_{p}=L^{p}\left( R^{n};E\right) \text{, }Y^{s,p}=H^{s,p}\left(
R^{n};E\right) ,\text{ }Y_{1}^{s,p}\left( A^{\gamma }\right) = 
\]

\[
H^{s,p}\left( R^{n};E\left( A^{\gamma }\right) \right) \cap L^{1}\left(
R^{n};E\left( A^{\gamma }\right) \right) ,\text{ for }0\leq \gamma \leq 1,%
\text{ } 
\]%
\[
Y_{\infty }^{s,p}=H^{s,p}\left( R^{n};E\right) \cap L^{\infty }\left(
R^{n};E\right) ,\text{ }Y_{\infty }^{s,p}\left( A\right) =H^{s,p}\cap
L^{\infty }\left( R^{n};E\left( A\right) \right) . 
\]

First we need the following lemmas

\textbf{Lemma 2.1. }Suppose the Condition 2.1 hold. Moreover, $\varphi ,$ $%
\psi $ $\in Y_{1}^{s,p}\left( A\right) $. Then problem $\left( 2.1\right)
-\left( 2.2\right) $ has a unique generalized solution.

\textbf{Proof. }By using of the Fourier transform we get from $(2.1)-(2.2):$%
\[
\hat{u}_{tt}\left( t,\xi \right) +A_{\xi }\hat{u}\left( t,\xi \right) =\hat{g%
}\left( t,\xi \right) ,\text{ } 
\]%
\begin{equation}
\hat{u}\left( 0,\xi \right) =\hat{\varphi}\left( \xi \right)
+\dint\limits_{0}^{T}\alpha \left( \sigma \right) \hat{u}\left( \xi ,\sigma
\right) d\sigma ,\text{ }  \tag{2.3}
\end{equation}

\[
\hat{u}_{t}\left( 0,\xi \right) =\hat{\psi}\left( \xi \right)
+\dint\limits_{0}^{T}\beta \left( \sigma \right) \hat{u}\left( \sigma ,\xi
\right) d\sigma ,\text{ }\xi \in R^{n},\text{ }t\in \left( 0,T\right) , 
\]%
where $\hat{u}\left( \xi ,t\right) $ is a Fourier transform of $u\left(
x,t\right) $ with respect to $x,$ where%
\[
A_{\xi }=A+\left\vert \xi \right\vert ^{2}\text{, }\xi \in R^{n}. 
\]%
Consider the problem%
\begin{equation}
\hat{u}_{tt}\left( t,\xi \right) +A_{\xi }\hat{u}\left( t,\xi \right) =\hat{g%
}\left( t,\xi \right) ,\text{ }  \tag{2.4}
\end{equation}%
\[
\hat{u}\left( \xi ,0\right) =u_{0}\left( \xi \right) ,\text{ }\hat{u}%
_{t}\left( \xi ,0\right) =u_{1}\left( \xi \right) ,\text{ }\xi \in R^{n},%
\text{ }t\in \left[ 0,T\right] ,\text{ } 
\]%
where $u_{0}\left( \xi \right) \in D\left( A\right) $ and $u_{1}\left( \xi
\right) \in D\left( A^{\frac{1}{2}}\right) $ for $\xi \in R^{n}.$ By virtue
of $\left[ \text{11, \S 11.2, 11.4}\right] $ we obtain that $A_{\xi }$ is a
generator of a strongly continuous cosine operator function\ and problem $%
(2.4)$ has a unique solution for all $\xi \in R^{n},$ moreover, the solution
can be written as%
\begin{equation}
\hat{u}\left( \xi ,t\right) =C\left( \xi ,t,A\right) u_{0}\left( \xi \right)
+S\left( \xi ,t,A\right) u_{1}\left( \xi \right) +  \tag{2.5}
\end{equation}%
\[
\dint\limits_{0}^{t}S\left( \xi ,t-\tau ,A\right) \hat{g}\left( \xi ,\tau
\right) d\tau ,\text{ }t\in \left( 0,T\right) , 
\]%
where $C\left( t,\xi ,A\right) $ is a cosine and $S\left( t,\xi ,A\right) $
is a sine operator-functions (see e.g. $\left[ 11\right] $) with generator
of $A_{\xi }$, i.e.%
\[
C\left( t,\xi ,A\right) =\frac{1}{2}\left( e^{itA_{\xi }^{\frac{1}{2}%
}}+e^{-itA_{\xi }^{\frac{1}{2}}}\right) \text{, }S\left( t,\xi ,A\right) =%
\frac{1}{2i}A_{\xi }^{-\frac{1}{2}}\left( e^{itA_{\xi }^{\frac{1}{2}%
}}-e^{-itA_{\xi }^{\frac{1}{2}}}\right) . 
\]%
Using formula $\left( 2.5\right) $ and nonlocal boundary condition 
\[
u_{0}\left( \xi \right) =\hat{\varphi}\left( \xi \right)
+\dint\limits_{0}^{T}\alpha \left( \sigma \right) \hat{u}\left( \xi ,\sigma
\right) d\sigma ,\text{ } 
\]%
we get 
\[
u_{0}\left( \xi \right) =\hat{\varphi}\left( \xi \right)
+\dint\limits_{0}^{T}\alpha \left( \sigma \right) \left[ C\left( \xi ,\sigma
,A\right) u_{0}\left( \xi \right) +S\left( \xi ,\sigma ,A\right) u_{1}\left(
\xi \right) \right] d\sigma + 
\]

\[
\dint\limits_{0}^{T}\dint\limits_{0}^{\sigma }S\left( \xi ,\sigma -\tau
,A\right) \hat{g}\left( \xi ,\tau \right) d\tau d\sigma ,\text{ }\tau \in
\left( 0,T\right) . 
\]%
Then, 
\[
\left[ I-\dint\limits_{0}^{T}\alpha \left( \sigma \right) C\left( \xi
,\sigma ,A\right) d\sigma \right] u_{0}\left( \xi \right) -\left[
\dint\limits_{0}^{T}\alpha \left( \sigma \right) S\left( \xi ,\sigma
,A\right) d\sigma \right] u_{1}\left( \xi \right) = 
\]%
\[
\dint\limits_{0}^{T}\dint\limits_{0}^{\sigma }\alpha \left( \sigma \right)
S\left( \xi ,\sigma -\tau ,A\right) \hat{g}\left( \xi ,\tau \right) d\tau
d\sigma +\hat{\varphi}\left( \xi \right) . 
\]%
Differentiating both sides of formula $\left( 2.5\right) $ we obtain%
\[
\hat{u}_{t}\left( \xi ,t\right) =-AS\left( \xi ,t,A\right) u_{0}\left( \xi
\right) +C\left( \xi ,t,A\right) u_{1}\left( \xi \right) + 
\]%
\[
\dint\limits_{0}^{t}C\left( \xi ,t-\tau ,A\right) \hat{g}\left( \xi ,\tau
\right) d\tau ,\text{ }t\in \left( 0,\infty \right) . 
\]%
Using this formula and integral condition%
\[
u_{1}\left( \xi \right) =\hat{\psi}\left( \xi \right)
+\dint\limits_{0}^{T}\beta \left( \sigma \right) \hat{u}_{t}\left( \xi
,\sigma \right) d\sigma 
\]%
we obtain%
\[
u_{1}\left( \xi \right) =\hat{\psi}\left( \xi \right)
+\dint\limits_{0}^{T}\beta \left( \sigma \right) \left[ -AS\left( \xi
,\sigma ,A\right) u_{0}\left( \xi \right) +C\left( \xi ,\sigma ,A\right)
u_{1}\left( \xi \right) \right] d\sigma + 
\]

\begin{equation}
\dint\limits_{0}^{T}\dint\limits_{0}^{\sigma }C\left( \xi ,\sigma -\tau
,A\right) \hat{g}\left( \xi ,\tau \right) d\tau d\sigma .  \tag{2.6}
\end{equation}%
Thus, 
\[
\dint\limits_{0}^{T}\beta \left( \sigma \right) AS\left( \xi ,\sigma
,A\right) d\sigma u_{0}\left( \xi \right) +\left[ I-\dint\limits_{0}^{T}%
\beta \left( \sigma \right) C\left( \xi ,\sigma ,A\right) d\sigma \right]
u_{1}\left( \xi \right) = 
\]%
\begin{equation}
\dint\limits_{0}^{T}\dint\limits_{0}^{\sigma }\beta \left( \sigma \right)
C\left( \xi ,\sigma -\tau ,A\right) \hat{g}\left( \xi ,\tau \right) d\tau
d\sigma +\hat{\psi}\left( \xi \right) .  \tag{2.7}
\end{equation}%
Now, we consider the system of equations $\left( 2.6\right) $ and $\left(
2.7\right) $ in $u_{0}\left( \xi \right) $ and $u_{1}\left( \xi \right) $.
The determinant of this system is 
\[
D\left( \xi \right) =\left\vert 
\begin{array}{cc}
\alpha _{11}\left( \xi \right) & \alpha _{12}\left( \xi \right) \\ 
\alpha _{21}\left( \xi \right) & \alpha _{22}\left( \xi \right)%
\end{array}%
\right\vert , 
\]%
where 
\[
\alpha _{11}\left( \xi \right) =I-\dint\limits_{0}^{T}\alpha \left( \sigma
\right) C\left( \xi ,\sigma ,A\right) d\sigma ,\text{ }\alpha _{12}\left(
\xi \right) =-\dint\limits_{0}^{T}\alpha \left( \sigma \right) S\left( \xi
,\sigma ,A\right) d\sigma , 
\]%
\[
\alpha _{21}\left( \xi \right) =\dint\limits_{0}^{T}\beta \left( \sigma
\right) AS\left( \xi ,\sigma ,A\right) d\sigma ,\text{ }\alpha _{22}\left(
\xi \right) =I-\dint\limits_{0}^{T}\beta \left( \sigma \right) S\left( \xi
,\sigma ,A\right) d\sigma . 
\]%
Then by using the properties%
\[
\left[ C\left( \sigma ,A\right) C\left( \tau ,A\right) +AS\left( \sigma
,A\right) S\left( \tau ,A\right) \right] =I 
\]%
of sine and cosine operator function $\left[ \text{11, \S 11.2, 11.4}\right] 
$ we obtain 
\[
D\left( \xi \right) =I-\dint\limits_{0}^{T}\left[ \alpha \left( \sigma
\right) +\beta \left( \sigma \right) \right] C\left( \sigma \right) d\sigma
+ 
\]%
\[
\dint\limits_{0}^{T}\dint\limits_{0}^{T}\alpha \left( \sigma \right) \beta
\left( \tau \right) \left[ C\left( \xi ,\sigma ,A\right) C\left( \xi ,\tau
,A\right) +AS\left( \xi ,\sigma ,A\right) S\left( \xi ,\tau A\right) \right]
d\sigma d\tau = 
\]

\[
I-\dint\limits_{0}^{T}\left[ \alpha \left( \sigma \right) +\beta \left(
\sigma \right) \right] C\left( \sigma \right) d\sigma
+\dint\limits_{0}^{T}\dint\limits_{0}^{T}\alpha \left( \sigma \right) \beta
\left( \tau \right) d\sigma d\tau =O\left( \xi \right) .
\]%
Solving the system $\left( 2.6\right) -\left( 2.7\right) $, we get%
\begin{equation}
u_{0}\left( \xi \right) =O^{-1}\left( \xi \right) \left\{ \left[
I-\dint\limits_{0}^{T}\beta \left( \sigma \right) C\left( \xi ,\sigma
,A\right) d\sigma \right] f_{1}\right. +\left. \dint\limits_{0}^{T}\alpha
\left( \sigma \right) S\left( \xi ,\sigma ,A\right) d\sigma f_{2}\right\} , 
\tag{2.9}
\end{equation}
\[
u_{1}\left( \xi \right) =O^{-1}\left( \xi \right) \left\{ \left[
I-\dint\limits_{0}^{T}\alpha \left( \sigma \right) C\left( \xi ,\sigma
,A\right) d\sigma \right] f_{2}\right. -\left. \dint\limits_{0}^{T}\left[
\beta \left( \sigma \right) A_{\xi }S\left( \xi ,\sigma ,A\right) d\sigma %
\right] f_{1}\right\} ,
\]%
where 
\[
f_{1}=\dint\limits_{0}^{T}\dint\limits_{0}^{\sigma }\alpha \left( \sigma
\right) S\left( \xi ,\sigma -\tau ,A\right) \hat{g}\left( \xi ,\tau \right)
d\tau d\sigma +\hat{\varphi}\left( \xi \right) ,
\]%
\begin{equation}
f_{2}=\dint\limits_{0}^{T}\dint\limits_{0}^{\sigma }\beta \left( \sigma
\right) C\left( \xi ,\sigma -\tau ,A\right) \hat{g}\left( \tau ,\xi \right)
d\tau d\sigma +\hat{\psi}\left( \xi \right) .  \tag{2.10}
\end{equation}%
From $\left( 2.5\right) ,$ $\left( 2.8\right) $ and $\left( 2.9\right) $ we
get that, the solution of the problem $\left( 2.4\right) $ can be expressed
as 
\[
\hat{u}\left( t,\xi \right) =O^{-1}\left( \xi \right) \left\{ C\left( \xi
,t,A\right) \left[ \left( I-\dint\limits_{0}^{T}\beta \left( \sigma \right)
C\left( \xi ,\sigma ,A\right) d\sigma \right) f_{1}\right. \right. +
\]%
\[
\left. \dint\limits_{0}^{T}\alpha \left( \sigma \right) S\left( \xi ,\sigma
,A\right) d\sigma f_{2}\right] +S\left( t,\xi ,A\right) \left[ \left(
I-\dint\limits_{0}^{T}\alpha \left( \sigma \right) C\left( \xi ,\sigma
,A\right) d\sigma \right) f_{2}\right. -
\]%
\begin{equation}
\left. \left. \dint\limits_{0}^{T}\beta \left( \sigma \right) A_{\xi
}S\left( \xi ,\sigma ,A\right) d\sigma f_{1}\right] \right\}
+\dint\limits_{0}^{t}S\left( t-\tau ,\xi ,A\right) \hat{g}\left( \tau ,\xi
\right) d\tau ,\text{ }t\in \left( 0,T\right) .  \tag{2.11}
\end{equation}%
Taking into account $\left( 2.10\right) $ we obtain from $\left( 2.11\right) 
$ that there is a generalized solution of $(2.1)-(2.2)$ given by 
\begin{equation}
u\left( x,t\right) =S_{1}\left( t,A\right) \varphi \left( x\right)
+S_{2}\left( t,A\right) \psi \left( x\right) +\Phi \left( x,t\right) , 
\tag{2.12}
\end{equation}%
where $S_{1}\left( t,A\right) $ and $S_{2}\left( t,A\right) $ are linear
operator functions in $E$ defined by 
\[
S_{1}\left( t,A\right) \varphi =\left( 2\pi \right) ^{-\frac{1}{n}%
}\dint\limits_{R^{n}}\left\{ e^{ix\xi }O^{-1}\left( \xi \right) \right. 
\text{ }
\]%
\[
\left[ C\left( t,\xi ,A\right) \left( I-\dint\limits_{0}^{T}\beta \left(
\sigma \right) C\left( \xi ,\sigma ,A\right) \right) -A_{\xi }S\left( \xi
,\sigma ,A\right) \right] d\sigma \left. \hat{\varphi}\left( \xi \right)
d\xi \right\} ,
\]

\begin{equation}
S_{2}\left( t,A\right) \psi =\left( 2\pi \right) ^{-\frac{1}{n}%
}\dint\limits_{R^{n}}\left\{ e^{ix\xi }O^{-1}\left( \xi \right) C\left(
t,\xi ,A\right) \right.  \tag{2.13}
\end{equation}%
\[
\dint\limits_{0}^{T}\left[ \alpha \left( \sigma \right) S\left( \xi ,\sigma
,A\right) +S\left( \xi ,\sigma ,A\right) \text{ }\left(
I-\dint\limits_{0}^{T}\alpha \left( \sigma \right) C\left( \xi ,\sigma
,A\right) \right) d\sigma \right] \left. \hat{\psi}\left( \xi \right)
\right\} d\xi 
\]%
\[
\left( I-\dint\limits_{0}^{T}\beta \left( \sigma \right) C\left( \xi ,\sigma
,A\right) -A_{\xi }S\left( \xi ,\sigma ,A\right) d\sigma \right) \left. \hat{%
\varphi}\left( \xi \right) \right\} d\xi , 
\]

\[
\Phi \left( x,t\right) =\left( 2\pi \right) ^{-\frac{1}{n}%
}\dint\limits_{R^{n}}O^{-1}\left( \xi \right) e^{ix\xi }\left\{
\dint\limits_{0}^{t}S\left( \xi ,t-\tau ,A\right) \hat{g}\left( \xi ,\tau
\right) d\tau \right. + 
\]

\[
\left[ C\left( \xi ,t,A\right) \left( I-\dint\limits_{0}^{T}\beta \left(
\sigma \right) C\left( \xi ,\sigma ,A\right) d\sigma \right) \right. + 
\]

\[
\left. S\left( \xi ,t,A\right) \dint\limits_{0}^{T}\beta \left( \sigma
\right) A_{\xi }S\left( \xi ,\sigma ,A\right) d\sigma \right] g_{1}\left(
\xi \right) +C\left( \xi ,t,A\right) \dint\limits_{0}^{T}\alpha \left(
\sigma \right) S\left( \xi ,\sigma ,A\right) d\sigma +
\]%
\[
S\left( \xi ,t,A\right) \left( I-\dint\limits_{0}^{T}\alpha \left( \sigma
\right) C\left( \xi ,\sigma ,A\right) d\sigma \right) \left. g_{2}\left( \xi
\right) \right\} d\xi ,
\]%
here 
\begin{equation}
g_{1}\left( \xi \right) =\dint\limits_{0}^{T}\dint\limits_{0}^{\sigma
}\alpha \left( \sigma \right) S\left( \xi ,\sigma -\tau ,A\right) \hat{g}%
\left( \xi ,\tau \right) d\tau d\sigma ,  \tag{2.14}
\end{equation}%
\[
g_{2}\left( \xi \right) =\dint\limits_{0}^{T}\dint\limits_{0}^{\sigma }\beta
\left( \sigma \right) C\left( \xi ,\sigma -\tau ,A\right) \hat{g}\left( \xi
,\tau \right) d\tau d\sigma .
\]

\textbf{Lemma 2.2. }Suppose the Condition 2.1 hold. Let $0\leq \gamma <\frac{%
1}{2}$ and $\varphi \in Y_{1}^{s,p}\left( A\right) $. Then the following
uniform estimate holds 
\begin{equation}
\left\Vert A^{\gamma }S_{1}\left( t,A\right) \varphi \right\Vert _{X_{\infty
}}\leq C_{1}\left[ \left\Vert A\varphi \right\Vert _{Y^{s,p}}+\left\Vert
A\varphi \right\Vert _{X_{1}}\right] .  \tag{2.15}
\end{equation}

\textbf{Proof. }Let $N\in \mathbb{N}$ and 
\[
\Pi _{N}=\left\{ \xi :\xi \in R^{n},\text{ }\left\vert \xi \right\vert \leq
N\right\} ,\text{ }\Pi _{N}^{\prime }=\left\{ \xi :\xi \in R^{n},\text{ }%
\left\vert \xi \right\vert \geq N\right\} . 
\]

It is clear to see that

\[
\left\Vert A^{\gamma }S_{1}\left( t,A\right) \varphi \right\Vert _{X_{\infty
}}=\left\Vert A^{\gamma }F^{-1}C\left( t,A\right) \hat{\varphi}\right\Vert
_{X_{\infty }}\leq 
\]

\[
\left\Vert F^{-1}A^{\gamma }C\left( \xi ,t,A\right) \hat{\varphi}\left( \xi
\right) \right\Vert _{L^{\infty }\left( \Pi _{N};E\right) }+\left\Vert
F^{-1}A^{\gamma }C\left( \xi ,t,A\right) \hat{\varphi}\left( \xi \right)
\right\Vert _{L^{\infty }\left( \Pi _{N}^{\prime };E\right) }\leq 
\]%
\begin{equation}
\left\Vert F^{-1}C\left( \xi ,t,A\right) A^{\gamma }\hat{\varphi}\left( \xi
\right) \right\Vert _{L^{\infty }\left( \Pi _{N};E\right) }+  \tag{2.16}
\end{equation}%
\[
\left\Vert F^{-1}\left( 1+\left\vert \xi \right\vert ^{2}\right) ^{-\frac{s}{%
2}}A^{-\left( 1-\gamma \right) }C\left( \xi ,t,A\right) \left( 1+\left\vert
\xi \right\vert ^{2}\right) ^{\frac{s}{2}}A\hat{\varphi}\left( \xi \right)
\right\Vert _{L^{\infty }\left( \Pi _{N}^{\prime };E\right) }. 
\]

In view of the Remark 1.2, properties of sectorial operator $A$ and by using
the Holder inequality we have%
\begin{equation}
\left\Vert F^{-1}A^{\gamma }C\left( \xi ,t,A\right) \hat{\varphi}\left( \xi
\right) \right\Vert _{L^{\infty }\left( \Pi _{N};E\right) }\leq
M_{1}\left\Vert A^{\gamma }\varphi \right\Vert _{X_{1}}.  \tag{2.17}
\end{equation}

By differentiating and in view of smoothness of $C\left( \xi ,t,A\right) $
in $\xi ,$ we have%
\[
\frac{\partial }{\partial \xi _{k}}\left[ \left( 1+\left\vert \xi
\right\vert ^{2}\right) ^{-\frac{s}{2}}A^{-\left( 1-\gamma \right) }C\left(
\xi ,t,A\right) \right] =-s\xi _{k}\left( 1+\left\vert \xi \right\vert
^{2}\right) ^{-\frac{s}{2}-1}A^{-\left( 1-\gamma \right) }C\left( \xi
,t,A\right) - 
\]%
\[
it\xi _{k}\left( 1+\left\vert \xi \right\vert ^{2}\right) ^{-\frac{s}{2}%
}A^{-\left( 1-\gamma \right) }S\left( \xi ,t,A\right) . 
\]

By differentiating again we obtain the same type operator-functons. By
applying $\left[ \text{11, Theorem 2.1}\right] $ and replacing $A$ by $%
A_{\xi }=A+\left\vert \xi \right\vert ^{2}$ for $\xi \in \Pi _{N}^{\prime }$
we obtain the unif\i rm estimate 
\[
\left\Vert C\left( \xi ,t,A\right) \right\Vert _{B\left( E\right) }\leq
M_{\xi }e^{\left\vert \omega \right\vert t},\text{ }\left\Vert S\left( \xi
,t,A\right) \right\Vert _{B\left( E\right) }\leq M_{\xi }e^{\left\vert
\omega \right\vert t}, 
\]%
where 
\[
M_{\xi }=M\left( 2+\lambda ^{2}\left\vert \left\vert \xi \right\vert
^{2}-\lambda ^{2}\right\vert ^{-1}\right) \leq 2M\text{ for all }\lambda \in
\left( 0,\infty \right) . 
\]

Moreover, by resolvent properties of sectorial operator $A$ we have 
\[
\left\Vert A_{\xi }^{-\left( 1-\gamma \right) }\right\Vert \leq C\left\vert
\xi \right\vert ^{-\left( 1-\gamma \right) }. 
\]

Hence, by above estimates of $C\left( \xi ,t,A\right) ,$ $S\left( \xi
,t,A\right) $ and resolvent properties of sectorial operator $A$ we obtain 
\begin{equation}
\left( 1+\left\vert \xi \right\vert ^{2}\right) ^{-\frac{s}{2}}A_{\xi
}^{-\left( 1-\gamma \right) }C\left( \xi ,t,A\right) \text{, }\left(
1+\left\vert \xi \right\vert ^{2}\right) ^{-\frac{s}{2}}A_{\xi }^{-\left(
1-\gamma \right) }C\left( \xi ,t,A\right) \in  \tag{2.18}
\end{equation}%
\[
W^{m,r}\left( R^{n};B\left( E\right) \right) \subset B_{r,1}^{n\left( \frac{1%
}{r}+\frac{1}{p}\right) }\left( R^{n};B\left( E\right) \right) 
\]%
for $m\geq s>n\left( \frac{1}{r}+\frac{1}{p}\right) $\ and $t\in \left[ 0,T%
\right] .$ By Proposition A$_{1}$ from $\left( 2.18\right) $ we get%
\begin{equation}
\left( 1+\left\vert \xi \right\vert ^{2}\right) ^{-\frac{s}{2}}A_{\xi
}^{-\left( 1-\gamma \right) }C\left( \xi ,t,A\right) \in M_{p}^{\infty
}\left( E\right) ,\text{ }  \tag{2.19}
\end{equation}%
\[
\left( 1+\left\vert \xi \right\vert ^{2}\right) ^{-\frac{s}{2}}A_{\xi
}^{-\left( 1-\gamma \right) }S\left( \xi ,t,A\right) \in M_{p}^{\infty
}\left( E\right) 
\]%
uniformly in $t\in \left[ 0,T\right] $, i.e. we have the following estimates%
\[
\left\Vert F^{-1}A^{\gamma }C\left( \xi ,t,A\right) \hat{\varphi}\left( \xi
\right) \right\Vert _{L^{\infty }\left( \Pi _{N}^{\prime };E\right) }\leq
M_{3}\left\Vert A\varphi \right\Vert _{Y^{s,p}}, 
\]%
\[
\left\Vert F^{-1}A^{\gamma }S\left( \xi ,t,A\right) \hat{\varphi}\left( \xi
\right) \right\Vert _{L^{\infty }\left( \Pi _{N}^{\prime };E\right) }\leq
M_{3}\left\Vert A\varphi \right\Vert _{Y^{s,p}}. 
\]

Then by Minkowski's inequality and from $\left( 2.16\right) $ and $\left(
2.17\right) $ we obtain $\left( 2.15\right) .$

\textbf{Lemma 2.3. }Suppose the Condition 2.1 hold. Let $0\leq \gamma <\frac{%
1}{2}$ and $\psi \in Y_{1}^{s,p}\left( A\right) .$ Then the uniform estimate
holds 
\begin{equation}
\left\Vert A^{\gamma }S_{2}\left( t,A\right) \psi \right\Vert _{X_{\infty
}}\leq M_{3}\left( \left\Vert A\psi \right\Vert _{Y^{s,p}}+\left\Vert A\psi
\right\Vert _{X_{1}}\right) .  \tag{2.20}
\end{equation}

\textbf{Proof. }It is clear that 
\[
\left\Vert A^{\gamma }S_{2}\left( t,A\right) \psi \right\Vert _{X_{\infty
}}=\left\Vert F^{-1}A^{\gamma }S\left( t,A\right) \hat{\psi}\right\Vert
_{X_{\infty }}\leq 
\]

\[
\left\Vert F^{-1}A^{\gamma }C\left( \xi ,t,A\right) \hat{\psi}\left( \xi
\right) \right\Vert _{L^{\infty }\left( \Pi _{N};E\right) }+\left\Vert
F^{-1}A^{\gamma }C\left( \xi ,t,A\right) \hat{\psi}\left( \xi \right)
\right\Vert _{L^{\infty }\left( \Pi _{N}^{\prime };E\right) }\leq 
\]%
\begin{equation}
\left\Vert F^{-1}S\left( \xi ,t,A\right) A^{-\left( 1-\gamma \right) }\hat{%
\psi}\left( \xi \right) \right\Vert _{L^{\infty }\left( \Pi _{N};E\right) }+
\tag{2.21}
\end{equation}%
\[
\left\Vert F^{-1}\left( 1+\left\vert \xi \right\vert ^{2}\right) ^{-\frac{s}{%
2}}A^{-\left( 1-\gamma \right) }S\left( \xi ,t,A\right) \left( 1+\left\vert
\xi \right\vert ^{2}\right) ^{\frac{s}{2}}A\hat{\psi}\left( \xi \right)
\right\Vert _{L^{\infty }\left( \Pi _{N}^{\prime };E\right) }. 
\]

\bigskip Then by $\left( 2.19\right) $ and by resolvent properties of
sectorial operator $A$ we obtain 
\begin{equation}
\left\Vert F^{-1}A^{\gamma }S\left( \xi ,t,A\right) \hat{\psi}\left( \xi
\right) \right\Vert _{L^{\infty }\left( \Pi _{N};E\right) }\leq
M_{4}\left\Vert A^{\gamma }\psi \right\Vert _{X_{1}},  \tag{2.22}
\end{equation}%
\[
\left\Vert F^{-1}A^{\gamma }S\left( \xi ,t,A\right) \hat{\psi}\left( \xi
\right) \right\Vert _{L^{\infty }\left( \Pi _{N}^{\prime };E\right) }\leq
M_{5}\left\Vert A\psi \right\Vert _{Y^{s,p}}. 
\]

Hence, from $\left( 2.21\right) $ and $\left( 2.22\right) $ we obtain $%
\left( 2.20\right) .$

\textbf{Lemma 2.4. }Suppose the Condition 2.1 hold. Let $0\leq \gamma <\frac{%
1}{2}$ and $g\left( .,t\right) \in Y_{1}^{s,p}$ for $t\in \left[ 0,T\right]
. $ Then we have the following uniform estimate 
\begin{equation}
\left\Vert \dint\limits_{0}^{t}A^{\gamma }S_{2}\left( x,t-\tau ,A\right)
g\left( x,\tau \right) d\tau \right\Vert _{X_{\infty }}\leq  \tag{2.23}
\end{equation}

\[
M_{6}\dint\limits_{0}^{t}\left( \left\Vert g\left( .,\tau \right)
\right\Vert _{Y^{s,p}}+\left\Vert g\left( .,\tau \right) \right\Vert
_{X_{1}}\right) d\tau . 
\]%
\textbf{Proof. }By reasoning as it has been made in the above, we get 
\begin{equation}
\left\Vert F^{-1}\dint\limits_{0}^{t}A^{\gamma }S\left( \xi ,t-\tau
,A\right) \hat{g}\left( \xi ,\tau \right) d\tau \right\Vert _{X_{\infty
}}\leq  \tag{2.24}
\end{equation}

\[
\left\Vert \dint\limits_{0}^{t}F^{-1}A^{\gamma }S\left( \xi ,t-\tau
,A\right) \hat{g}\left( \xi ,\tau \right) d\tau \right\Vert _{L^{\infty
}\left( \Pi _{N};E\right) }+ 
\]%
\[
\left\Vert \dint\limits_{0}^{t}F^{-1}\left( 1+\left\vert \xi \right\vert
^{2}\right) ^{-\frac{s}{2}}A^{\gamma }S\left( \xi ,t-\tau ,A\right) \left(
1+\left\vert \xi \right\vert ^{2}\right) ^{\frac{s}{2}}\hat{g}\left( \xi
,\tau \right) d\tau \right\Vert _{L^{\infty }\left( \Pi _{N}^{\prime
};E\right) }. 
\]

In view of smoothness $C\left( \xi ,t,A\right) $ with respect to $\xi \in
R^{n},$ resolvent properties of sectorial operator $A$ and by Remark 1.2 we
get 
\[
\left( 1+\left\vert \xi \right\vert ^{2}\right) ^{-\frac{s}{2}%
}\dint\limits_{0}^{t}A^{\gamma }S\left( \xi ,t-\tau ,A\right) d\tau \in
B_{q,1}^{n\left( \frac{1}{r}+\frac{1}{p}\right) }\left( R^{n};B\left(
E\right) \right) 
\]%
for $s>n\left( \frac{1}{r}+\frac{1}{p}\right) $\ for all $t\in \left[ 0,T%
\right] .$

\bigskip Then by Proposition A$_{1}$ we get that, the operator-valued
functions 
\[
\left( 1+\left\vert \xi \right\vert ^{2}\right) ^{-\frac{s}{2}%
}\dint\limits_{0}^{t}A^{\gamma }S\left( \xi ,t-\tau ,A\right) d\tau \in
M_{p}^{\infty }\left( E\right) 
\]
uniformly in $t\in \left[ 0,T\right] .$ Hence, $\left( 2.24\right) $ implies 
$\left( 2.23\right) .$

\textbf{Theorem 2.1. }Assume the Condition 2.1 hold and $0\leq \gamma <\frac{%
1}{2}$. Moreover, $\varphi \in Y_{1}^{s,p}\left( A\right) ,$ $\psi \in
Y_{1}^{s,p}\left( A\right) $ and $g\left( .,t\right) \in Y_{1}^{s,p}$ for $%
t\in \left[ 0,T\right] .$ Then problem $\left( 2.1\right) -\left( 2.2\right) 
$\ has a unique solution $u\left( x,t\right) \in C^{2}\left( \left[ 0,T%
\right] ;Y_{\infty }^{s,p}\left( A\right) \right) $ and the following
estimate holds 
\begin{equation}
\left\Vert A^{\gamma }u\right\Vert _{X_{\infty }}+\left\Vert A^{\gamma
}u_{t}\right\Vert _{X_{\infty }}\leq C\left\{ \left\Vert A\varphi
\right\Vert _{Y^{s,p}}+\left\Vert A\varphi \right\Vert _{X_{1}}\right. + 
\tag{2.25}
\end{equation}

\[
\left\Vert A\psi \right\Vert _{Y^{s,p}}+\left\Vert A\psi \right\Vert
_{X_{1}}+\left. \dint\limits_{0}^{t}\left( \left\Vert \Delta g\left( .,\tau
\right) \right\Vert _{Y^{s,p}}+\left\Vert \Delta g\left( .,\tau \right)
\right\Vert _{X_{1}}\right) d\tau \right\} 
\]%
uniformly in $t\in \left[ 0,T\right] .$

\textbf{Proof. }From Lemma 2.1 we obtain that the problem $\left( 2.1\right)
-\left( 2.2\right) $\ has a solution $u$. From $\left( 2.12\right) $ and
from the estimates $\left( 2.15\right) ,$ $\left( 2.20\right) ,$ $\left(
2.23\right) $ for $0\leq \gamma <\frac{1}{2}$ we get the estimate 
\[
\left\Vert A^{\gamma }u\left( .,t\right) \right\Vert _{X_{\infty }}\leq
C\left\{ \left\Vert A\varphi \right\Vert _{Y^{s,p}}+\left\Vert A\varphi
\right\Vert _{X_{1}}+\left\Vert A\psi \right\Vert _{Y^{s,p}}+\left\Vert
A\psi \right\Vert _{X_{1}}\right. + 
\]

\begin{equation}
+\left. \dint\limits_{0}^{t}\left( \left\Vert \Delta g\left( .,\tau \right)
\right\Vert _{Y^{s,p}}+\left\Vert \Delta g\left( .,\tau \right) \right\Vert
_{X_{1}}\right) d\tau \right\} .  \tag{2.26}
\end{equation}%
By differentiating, in view of $\left( 2.5\right) $ we get from $\left(
2.12\right) $ the estimate of type $\left( 2.26\right) $ for $u_{t},$ when $%
0\leq \gamma <\frac{1}{2}.$ From this we obtain the estimate $\left(
2.25\right) .$

\textbf{Theorem 2.2. }Assume the Condition 2.1 hold and $0\leq \gamma <\frac{%
1}{2}$. Moreover, $\varphi \in Y_{1}^{s,p}\left( A\right) $, $\psi \in
Y_{1}^{s,p}\left( A\right) $ and $g\left( .,t\right) \in Y^{s,p}$ for $t\in
\left( 0,\infty \right) .$ Then the problem $\left( 2.1\right) -\left(
2.2\right) $\ has a unique solution $u\left( x,t\right) \in C^{2}\left( %
\left[ 0,T\right] ;Y^{s,p}\left( A\right) \right) $ and the following
uniform estimate holds%
\begin{equation}
\left( \left\Vert A^{\gamma }u\left( .,t\right) \right\Vert
_{Y^{s,p}}+\left\Vert A^{\gamma }u_{t}\left( .,t\right) \right\Vert
_{Y^{s,p}}\right) \leq  \tag{2.27}
\end{equation}

\[
C\left( \left\Vert A\varphi \right\Vert _{Y^{s,p}}+\left\Vert A\psi
\right\Vert _{Y^{s,p}}+\dint\limits_{0}^{t}\left\Vert \Delta g\left( \tau
,.\right) \right\Vert _{Y^{s,p}}d\tau \right) .
\]%
\textbf{Proof. }From $\left( 2.5\right) $ we have the following estimate 
\begin{equation}
\left( \left\Vert F^{-1}\left( 1+\left\vert \xi \right\vert ^{2}\right) ^{%
\frac{s}{2}}\hat{u}\right\Vert _{X_{p}}+\left\Vert F^{-1}\left( 1+\left\vert
\xi \right\vert ^{2}\right) ^{\frac{s}{2}}\hat{u}_{t}\right\Vert
_{X_{p}}\right) \leq   \tag{2.28}
\end{equation}

\[
C\left\{ \left\Vert F^{-1}\left( 1+\left\vert \xi \right\vert ^{2}\right) ^{%
\frac{s}{2}}C\left( \xi ,t,A\right) \hat{\varphi}\right\Vert _{X_{p}}\right.
+\left\Vert F^{-1}\left( 1+\left\vert \xi \right\vert ^{2}\right) ^{\frac{s}{%
2}}S\left( \xi ,t,A\right) \hat{\psi}\right\Vert _{X_{p}}+
\]

\[
\left. \dint\limits_{0}^{t}\left\Vert F^{-1}\left( 1+\left\vert \xi
\right\vert ^{2}\right) ^{\frac{s}{2}}S\left( t-\tau ,\xi ,A\right) \hat{g}%
\left( .,\tau \right) \right\Vert _{X_{p}}d\tau \right\} .
\]

\bigskip By above estimates of $C\left( \xi ,t,A\right) ,$ $S\left( \xi
,t,A\right) $ and resolvent properties of sectorial operator $A$ we obtain 
\[
\left( 1+\left\vert \xi \right\vert ^{2}\right) ^{-\frac{s}{2}}A_{\xi
}^{-\left( 1-\gamma \right) }C\left( \xi ,t,A\right) \text{, }\left(
1+\left\vert \xi \right\vert ^{2}\right) ^{-\frac{s}{2}}A_{\xi }^{-\left(
1-\gamma \right) }C\left( \xi ,t,A\right) \in B_{r,1}^{\frac{n}{r}}\left(
R^{n};B\left( E\right) \right) .
\]

\bigskip Then by Proposition A$_{1}$ from $\left( 2.18\right) $ we get%
\begin{equation}
\left( 1+\left\vert \xi \right\vert ^{2}\right) ^{-\frac{s}{2}}A_{\xi
}^{-\left( 1-\gamma \right) }C\left( \xi ,t,A\right) \in M_{p}^{p}\left(
E\right) ,\text{ }  \tag{2.29}
\end{equation}%
\[
\left( 1+\left\vert \xi \right\vert ^{2}\right) ^{-\frac{s}{2}}A_{\xi
}^{-\left( 1-\gamma \right) }S\left( \xi ,t,A\right) \in M_{p}^{p}\left(
E\right) 
\]%
for all $t$. So, the relations $\left( 2.29\right) $ by using the
Minkowski's inequality for integrals implies $\left( 2.27\right) .$

\begin{center}
\textbf{3. Initial value problem for nonlinear equation}
\end{center}

In this section, we will show the local existence and uniqueness of solution
to Cauchy problem $(1.1)-(1.2).$

For the study of the nonlinear problem $\left( 1.1\right) -\left( 1.2\right) 
$ we need the following lemma from $\left[ 31\right] .$

\textbf{Lemma 3.1} (Abstract Nirenberg's inequality). Let $E$ be a $UMD$
space. Assume that $u\in L^{p}\left( R^{n};E\right) $, $D^{m}u$ $\in
L^{q}\left( R^{n};E\right) $, $p,q\in \left( 1,\infty \right) $. Then for $i$
with $0\leq i\leq m,$ $m>\frac{n}{q}$ we have 
\[
\left\Vert D^{i}u\right\Vert _{r}\leq C\left\Vert u\right\Vert _{p}^{1-\mu
}\dsum\limits_{k=1}^{n}\left\Vert D_{k}^{m}u\right\Vert _{q}^{\mu },
\]%
where%
\[
\frac{1}{r}=\frac{i}{m}+\mu \left( \frac{1}{q}-\frac{m}{n}\right) +\left(
1-\mu \right) \frac{1}{p},\text{ }\frac{i}{m}\leq \mu \leq 1.
\]

Note that, for $E=\mathbb{C}$ the lemma considered by L. Nirenberg $\left[ 31%
\right] .$

Using the chain rule of the composite function, from Lemma 3.1 we obtain the
following result

\textbf{Lemma 3.2. }Let $E$ be an $UMD$ space. Assume that $u\in $ $%
W^{s,p}\left( R^{n};E\right) \cap L^{\infty }\left( R^{n};E\right) $, and $%
F\left( u\right) $ possesses continuous derivatives up to order $s\geq 1$.
Then $F\left( u\right) -F\left( 0\right) \in W^{m,p}\left( \Omega ;E\right) $
and 
\[
\left\Vert F\left( u\right) -F\left( 0\right) \right\Vert _{p}\leq
\left\Vert F^{^{\left( 1\right) }}\left( u\right) \right\Vert _{\infty
}\left\Vert u\right\Vert _{p}, 
\]

\[
\left\Vert D^{k}F\left( u\right) \right\Vert _{p}\leq
C_{0}\dsum\limits_{j=1}^{k}\left\Vert F^{\left( j\right) }\left( u\right)
\right\Vert _{\infty }\left\Vert u\right\Vert _{\infty }^{j-1}\left\Vert
D^{k}u\right\Vert _{p}\text{, }1\leq k\leq s, 
\]%
where $C_{0}$ $\geq 1$ is a constant and $k$ is an integer number.

For $E=\mathbb{C}$ the lemma coincide with the corresponding inequality in $%
\left[ 22\right] .$ Let%
\[
\text{ }X=L^{p}\left( R^{n};E\right) ,\text{ }Y=W^{s,p}\left( R^{n};E\left(
A\right) ,E\right) ,\text{ }E_{0}=\left( X,Y\right) _{\frac{1}{2p},p}. 
\]

\bigskip \textbf{Remark 3.1. }By using J. Lions-I. Petree result (see e.g $%
\left[ \text{27}\right] $ or $\left[ \text{40, \S\ 1.8.}\right] $) we obtain
that the map $u\rightarrow u\left( t_{0}\right) $, $t_{0}\in \left[ 0,T%
\right] $ is continuous from $W^{s,p}\left( 0,T;X,Y\right) $ onto $E_{0}$
and there is a constant $C_{1}$ such that 
\[
\left\Vert u\left( t_{0}\right) \right\Vert _{E_{0}}\leq C_{1}\left\Vert
u\right\Vert _{W^{s,p}\left( 0,T;X,Y\right) },\text{ }1\leq p\leq \infty 
\text{.} 
\]%
Let we define the space $Y\left( T\right) =C\left( \left[ 0,T\right]
;Y_{\infty }^{s,p}\left( A\right) \right) $ equipped with the norm defined by%
\[
\left\Vert u\right\Vert _{Y\left( T\right) }=\max\limits_{t\in \left[ 0,T%
\right] }\left\Vert u\right\Vert _{Y^{s,p}}+\max\limits_{t\in \left[ 0,T%
\right] }\left\Vert u\right\Vert _{L^{\infty }\left( R^{n};E\left( A\right)
\right) },\text{ }u\in Y\left( T\right) . 
\]%
It is easy to see that $Y\left( T\right) $ is a Banach space. For $\varphi
\in Y_{\infty }^{s,p}\left( A\right) $ and $\psi \in Y_{\infty }^{s,p}\left(
A^{\frac{1}{2}}\right) $ let 
\[
M=\left\Vert A\varphi \right\Vert _{Y^{s,p}}+\left\Vert A\varphi \right\Vert
_{X_{\infty }}+\left\Vert A\psi \right\Vert _{Y^{s,p}}+\left\Vert A\psi
\right\Vert _{X_{\infty }}. 
\]

\textbf{Definition 3.1. }For any $T>0$ if $\varphi \in Y_{\infty
}^{s,p}\left( A\right) $, $\psi \in Y_{\infty }^{s,p}\left( A\right) $ and $%
u $ $\in C\left( \left[ 0,T\right] ;Y_{\infty }^{s,p}\left( A\right) \right) 
$ satisfies the equation $(1.1)-(1.2)$ then $u\left( x,t\right) $ is called
the continuous solution\ or the strong solution of the problem $(1.1)-(1.2).$
If $T<\infty $, then $u\left( x,t\right) $ is called the local strong
solution of the problem $(1.1)-(1.2).$ If $T=\infty $, then $u\left(
x,t\right) $ is called the global strong solution of the problem $%
(1.1)-(1.2) $.

\textbf{Condition 3.1. }Assume:

(1) $E$ is a UMD space and $A\in \sigma \left( C_{0},\omega ,E\right) $;

(2) $\varphi \in Y_{\infty }^{s,p}\left( A\right) $, $\psi \in Y_{\infty
}^{s,p}\left( A\right) $ and $1\leq p<\infty $ for $s>n\left( \frac{1}{r}+%
\frac{1}{p}\right) $;

(3) the function $u\rightarrow $ $F\left( u\right) $: $R^{n}\times \left[ 0,T%
\right] \times E_{0}\rightarrow E$ is a measurable in $\left( x,t\right) \in
R^{n}\times \left[ 0,T\right] $ for $u\in E_{0};$

(4) $F\left( x,t,.,.\right) $ is continuous in $u\in E_{0}$ for $x\in R^{n},$
$t\in \left[ 0,T\right] $ and $f\left( u\right) \in C^{\left( 1\right)
}\left( E_{0};E\right) .$

\bigskip Main aim of this section is to prove the following result:

\bigskip \textbf{Theorem 3.1. }Assume the Condition 3.1 are satisfied. Then
problem $\left( 1.1\right) -\left( 2.2\right) $ has a unique local strange
solution $u\in C^{\left( 2\right) }\left( \left[ 0\right. ,\left.
T_{0}\right) ;Y_{\infty }^{s,p}\left( A\right) \right) ,$ where $T_{0}$ is a
maximal time interval that is appropriately small relative to $M$. Moreover,
if

\begin{equation}
\sup_{t\in \left[ 0\right. ,\left. T_{0}\right) }\left( \left\Vert
u\right\Vert _{Y^{s,p}}+\left\Vert Au\right\Vert _{X_{\infty }}+\left\Vert
Au_{t}\right\Vert _{Y^{s,p}}+\left\Vert Au_{t}\right\Vert _{X_{\infty
}}\right) <\infty  \tag{3.0}
\end{equation}%
then $T_{0}=\infty .$

\textbf{Proof. }First, we are going to prove the existence and the
uniqueness of the local continuous solution of the problem $(1.1)-\left(
1.2\right) $ by contraction mapping principle. Suppose that $u\in C^{\left(
2\right) }\left( \left[ 0,T\right] ;Y_{\infty }^{s,p}\left( A\right) \right) 
$ is a strong solution of the problem $(1.1)-(1.2)$. Consider a map $G$ on $%
Y\left( T\right) $ such that $G(u)$ is the solution of the Cauchy problem%
\begin{equation}
G_{tt}\left( u\right) -\Delta G_{tt}\left( u\right) +AG\left( u\right)
=F\left( G\left( u\right) \right) ,\text{ }x\in R^{n},\text{ }t\in \left(
0,T\right) ,  \tag{3.1}
\end{equation}%
\[
G\left( u\right) \left( x,0\right) =\varphi \left( x\right)
+\dint\limits_{0}^{T}\alpha \left( \sigma \right) G\left( u\right) \left(
x,\sigma \right) d\sigma ,\text{ } 
\]

\[
G\left( u\right) _{t}\left( x,0\right) =\psi \left( x\right)
+\dint\limits_{0}^{T}\beta \left( \sigma \right) G\left( u\right) _{t}\left(
x,\sigma \right) d\sigma , 
\]%
From Lemma 3.2 we know that $F(u)\in $ $L^{p}\left( 0,T;Y_{\infty
}^{s,p}\right) $ for any $T>0$. Thus, in view of Remark 1.3 and by Theorem
2.1, problem $\left( 3.1\right) $ has a unique solution which can be written
as%
\begin{equation}
G\left( u\right) \left( t,x\right) =S_{1}\left( t,A\right) \varphi \left(
x\right) +S_{2}\left( t,A\right) \psi \left( x\right) +\tilde{\Phi}\left(
t,x\right) ,  \tag{3.2}
\end{equation}%
where $S_{1}\left( t,A\right) ,$ $S_{2}\left( t,A\right) $ are defined by $%
\left( 2.13.\right) $ and 
\begin{equation}
\tilde{\Phi}\left( t,x\right) =\left( 2\pi \right) ^{-\frac{1}{n}%
}\dint\limits_{R^{n}}O^{-1}\left( \xi \right) e^{ix\xi }\left\{
\dint\limits_{0}^{t}S\left( t-\tau ,\xi ,A\right) \hat{F}\left( u\right)
\left( \xi ,\tau \right) d\tau \right. +  \tag{3.3}
\end{equation}

\[
\left[ C\left( t,\xi ,A\right) \left( I-\dint\limits_{0}^{T}\beta \left(
\sigma \right) C\left( \sigma ,\xi ,A\right) d\sigma \right) \right. + 
\]

\[
\left. S\left( t,\xi ,A\right) \dint\limits_{0}^{T}\beta \left( \sigma
\right) A_{\xi }S\left( \sigma ,\xi ,A\right) d\sigma \right] g_{1}\left(
\xi \right) +\left[ C\left( t,\xi ,A\right) \dint\limits_{0}^{T}\alpha
\left( \sigma \right) S\left( \sigma ,\xi ,A\right) d\sigma +\right. 
\]%
\[
\left. \left. S\left( t,\xi ,A\right) \left( I-\dint\limits_{0}^{T}\alpha
\left( \sigma \right) C\left( \sigma ,\xi ,A\right) d\sigma \right)
g_{2}\left( \xi \right) \right] \right\} d\xi .
\]%
Here, 
\begin{equation}
\tilde{g}_{1}\left( \xi \right)
=\dint\limits_{0}^{T}\dint\limits_{0}^{\sigma }\alpha \left( \sigma \right)
S\left( \sigma -\tau ,\xi ,A\right) \hat{F}\left( u\right) \left( \tau ,\xi
\right) d\tau d\sigma ,  \tag{3.4}
\end{equation}%
\[
\tilde{g}_{2}\left( \xi \right)
=\dint\limits_{0}^{T}\dint\limits_{0}^{\sigma }\beta \left( \sigma \right)
C\left( \sigma -\tau ,\xi ,A\right) \hat{F}\left( u\right) \left( \tau ,\xi
\right) d\tau d\sigma .
\]%
For the sake of convenience, we assume that $F\left( 0\right) =0$.
Otherwise, we can replace $F\left( u\right) $ with $F\left( u\right)
-F\left( 0\right) $. Hence, from Lemma $3.2$ we have $F\left( u\right) $ $%
\in Y^{s,p}$ iff $F\in C^{\left( k\right) }\left( R;E\right) $ for $k\geq s>1
$. Consider the operator in $Y\left( T\right) $ defined as 
\begin{equation}
Gu=S_{1}\left( t,A\right) \varphi \left( x\right) +S_{2}\left( t,A\right)
\psi \left( x\right) +\tilde{\Phi}\left( t,x\right) .  \tag{3.5}
\end{equation}%
From Lemma 3.2 we get that the operator $G$ is well defined for $F\in $ $%
C\left( R;E\right) .$ Moreover, from Lemma 3.2 it is easy to see that the
map $G$ is well defined for $F\in C^{\left( k\right) }\left( X_{0};E\right) $
for $k\geq s>1$ and $k\in \mathbb{N}.$ We put 
\[
Q\left( M;T\right) =\left\{ u\mid u\in Y\left( T\right) \text{, }\left\Vert
u\right\Vert _{Y\left( T\right) }\leq M+1\right\} .
\]%
Let us prove that the map $G$ has a unique fixed point in $Q\left(
M;T\right) .$ For this aim, it is sufficient to show that the operator $G$
maps $Q\left( M;T\right) $ into $Q\left( M;T\right) $ and $G:$ $Q\left(
M;T\right) $ $\rightarrow $ $Q\left( M;T\right) $ is strictly contractive if 
$T$ is appropriately small relative to $M.$ Consider the function \ $\bar{f}%
\left( \xi \right) $: $\left[ 0,\right. $ $\left. \infty \right) \rightarrow %
\left[ 0,\right. $ $\left. \infty \right) $ defined by 
\[
\ \bar{f}\left( \xi \right) =\max\limits_{\left\vert x\right\vert \leq \xi
}\left\{ \left\Vert F^{\left( 1\right) }\left( x\right) \right\Vert _{E},...,%
\text{ }\left\Vert F^{\left( k\right) }\left( x\right) \right\Vert _{E}\text{
}\right\} ,\text{ }\xi \geq 0.
\]%
It is clear to see that the function $\bar{f}\left( \xi \right) $ is
continuous and nondecreasing on $\left[ 0,\right. $ $\left. \infty \right) .$
From Lemma 3.2 for $1\leq k\leq s$ we have

\[
\left\Vert F\left( u\right) \right\Vert _{Y^{k,p}}\leq \left\Vert
F^{^{\left( 1\right) }}\left( u\right) \right\Vert _{\infty }\left\Vert
u\right\Vert _{p}+C_{0}\dsum\limits_{j=1}^{k}\left\Vert F^{\left( j\right)
}\left( u\right) \right\Vert _{\infty }\left\Vert u\right\Vert _{\infty
}^{j-1}\left\Vert D^{k}u\right\Vert _{p}\leq 
\]%
\qquad\ 
\begin{equation}
2C_{0}\bar{f}\left( M+1\right) \left( M+1\right) \left\Vert u\right\Vert
_{Y^{s,p}}.  \tag{3.6}
\end{equation}%
By using the Theorem 2.1 we obtain from $\left( 3.5\right) $%
\begin{equation}
\left\Vert A^{\gamma }G\left( u\right) \right\Vert _{X_{\infty }}\leq
\left\Vert A\varphi \right\Vert _{X_{\infty }}+\left\Vert A^{\frac{1}{2}%
}\psi \right\Vert _{X_{\infty }}+\dint\limits_{0}^{t}\left\Vert \Delta
F\left( u\left( \tau \right) \right) \right\Vert _{X_{\infty }},  \tag{3.7}
\end{equation}%
\begin{equation}
\left\Vert A^{\gamma }G\left( u\right) \right\Vert _{Y^{2,p}}\leq \left\Vert
A\varphi \right\Vert _{Y^{2,p}}+\left\Vert A^{\frac{1}{2}}\psi \right\Vert
_{Y^{2,p}}+\dint\limits_{0}^{t}\left\Vert \Delta F\left( u\left( \tau
\right) \right) \right\Vert _{Y^{2,p}}d\tau .  \tag{3.8}
\end{equation}%
Thus, from $\left( 3.6\right) -\left( 3.8\right) $ and Lemma 3.2 we get 
\[
\left\Vert A^{\gamma }G\left( u\right) \right\Vert _{Y\left( T\right) }\leq
M+T\left( M+1\right) \left[ 1+2C_{0}\left( M+1\right) \bar{f}\left(
M+1\right) \right] .
\]%
If $T$ satisfies 
\begin{equation}
T\leq \left\{ \left( M+1\right) \left[ 1+2C_{0}\left( M+1\right) \bar{f}%
\left( M+1\right) \right] \right\} ^{-1}\left\Vert Gu\right\Vert _{Y\left(
T\right) }\leq M+1.  \tag{3.9}
\end{equation}%
then Therefore, if $\left( 3.9\right) $ holds, then $G$ maps $Q\left(
M;T\right) $ into $Q\left( M;T\right) .$ Now, we are going to prove that the
map $G$ is strictly contractive. Assume $T>0$ and $u_{1},$ $u_{2}\in $ $%
Q\left( M;T\right) $ given. We get%
\[
G\left( u_{1}\right) -G\left( u_{2}\right) =\left( 2\pi \right) ^{-\frac{1}{n%
}}\dint\limits_{R^{n}}\dint\limits_{0}^{t}\check{S}\left( t-\tau ,\xi
,A\right) \left[ \hat{F}\left( u_{1}\right) \left( \xi ,\tau \right) -\hat{F}%
\left( u_{2}\right) \left( \xi ,\tau \right) \right] d\tau d\xi .
\]

Using the mean value theorem, we obtain%
\[
\hat{F}\left( u_{1}\right) -\hat{F}\left( u_{2}\right) =\hat{F}^{\left(
1\right) }\left( u_{2}+\eta _{1}\left( u_{1}-u_{2}\right) \right) \left(
u_{1}-u_{2}\right) ,\text{ }
\]

\[
D_{\xi }\left[ \hat{F}\left( u_{1}\right) -\hat{F}\left( u_{2}\right) \right]
=\hat{F}^{\left( 2\right) }\left( u_{2}+\eta _{2}\left( u_{1}-u_{2}\right)
\right) \left( u_{1}-u_{2}\right) D_{\xi }u_{1}+\text{ } 
\]%
\[
\hat{F}^{\left( 1\right) }\left( u_{2}\right) \left( D_{\xi }u_{1}-D_{\xi
}u_{2}\right) , 
\]%
\[
D_{\xi }^{2}\left[ \hat{F}\left( u_{1}\right) -\hat{F}\left( u_{2}\right) %
\right] =\hat{F}^{\left( 3\right) }\left( u_{2}+\eta _{3}\left(
u_{1}-u_{2}\right) \right) \left( u_{1}-u_{2}\right) \left( D_{\xi
}u_{1}\right) ^{2}+\text{ } 
\]%
\[
\hat{F}^{\left( 2\right) }\left( u_{2}\right) \left( D_{\xi }u_{1}-D_{\xi
}u_{2}\right) \left( D_{\xi }u_{1}+D_{\xi }u_{2}\right) + 
\]%
\[
\hat{F}^{\left( 2\right) }\left( u_{2}+\eta _{4}\left( u_{1}-u_{2}\right)
\right) \left( u_{1}-u_{2}\right) D_{\xi }^{2}u_{1}+\hat{F}^{\left( 1\right)
}\left( u_{2}\right) \left( D_{\xi }^{2}u_{1}-D_{\xi }^{2}u_{2}\right) , 
\]%
where $0<\eta _{i}<1,$ $i=1,2,3,4.$ Thus using Holder's and Nirenberg's
inequality, we have%
\begin{equation}
\left\Vert \hat{F}\left( u_{1}\right) -\hat{F}\left( u_{2}\right)
\right\Vert _{X_{\infty }}\leq \bar{f}\left( M+1\right) \left\Vert
u_{1}-u_{2}\right\Vert _{X_{\infty }},  \tag{3.10}
\end{equation}%
\begin{equation}
\left\Vert \hat{F}\left( u_{1}\right) -\hat{F}\left( u_{2}\right)
\right\Vert _{X_{p}}\leq \bar{f}\left( M+1\right) \left\Vert
u_{1}-u_{2}\right\Vert _{X_{p}},  \tag{3.11}
\end{equation}%
\begin{equation}
\left\Vert D_{\xi }\left[ \hat{F}\left( u_{1}\right) -\hat{F}\left(
u_{2}\right) \right] \right\Vert _{X_{p}}\leq \left( M+1\right) \bar{f}%
\left( M+1\right) \left\Vert u_{1}-u_{2}\right\Vert _{X_{\infty }}+ 
\tag{3.12}
\end{equation}%
\[
\bar{f}\left( M+1\right) \left\Vert \hat{F}\left( u_{1}\right) -\hat{F}%
\left( u_{2}\right) \right\Vert _{X_{p}}, 
\]%
\[
\left\Vert D_{\xi }^{2}\left[ \hat{F}\left( u_{1}\right) -\hat{F}\left(
u_{2}\right) \right] \right\Vert _{X_{p}}\leq \left( M+1\right) \bar{f}%
\left( M+1\right) \left\Vert u_{1}-u_{2}\right\Vert _{X_{\infty }}\left\Vert
D_{\xi }^{2}u_{1}\right\Vert _{Y^{2,p}}^{2}+ 
\]%
\[
\bar{f}\left( M+1\right) \left\Vert D_{\xi }\left( u_{1}-u_{2}\right)
\right\Vert _{Y^{2,p}}\left\Vert D_{\xi }\left( u_{1}+u_{2}\right)
\right\Vert _{Y^{2,p}}+ 
\]%
\[
\bar{f}\left( M+1\right) \left\Vert u_{1}-u_{2}\right\Vert _{X_{\infty
}}\left\Vert D_{\xi }^{2}u_{1}\right\Vert _{X_{p}}+\bar{f}\left( M+1\right)
\left\Vert D_{\xi }\left( u_{1}-u_{2}\right) \right\Vert _{X_{p}}\leq 
\]%
\[
C^{2}\bar{f}\left( M+1\right) \left\Vert u_{1}-u_{2}\right\Vert _{X_{\infty
}}\left\Vert u_{1}\right\Vert _{X_{\infty }}\left\Vert D_{\xi
}^{2}u_{1}\right\Vert _{X_{p}}+C^{2}\bar{f}\left( M+1\right) 
\]%
\[
\left\Vert u_{1}-u_{2}\right\Vert _{X_{\infty }}^{\frac{1}{2}}\left\Vert
D_{\xi }^{2}\left( u_{1}-u_{2}\right) \right\Vert _{X_{p}}\left\Vert
u_{1}+u_{2}\right\Vert _{X_{\infty }}^{\frac{1}{2}}\left\Vert D_{\xi
}^{2}\left( u_{1}+u_{2}\right) \right\Vert _{X_{p}}+ 
\]%
\[
\left( M+1\right) \bar{f}\left( M+1\right) \left\Vert u_{1}-u_{2}\right\Vert
_{X_{\infty }}+\bar{f}\left( M+1\right) \left\Vert D_{\xi }^{2}\left(
u_{1}-u_{2}\right) \right\Vert _{X_{p}}\leq 3C^{2}\left( M+1\right) ^{2} 
\]%
\[
\bar{f}\left( M+1\right) \left\Vert u_{1}-u_{2}\right\Vert _{X_{\infty
}}+2C^{2}\left( M+1\right) \bar{f}\left( M+1\right) \left\Vert D_{\xi
}^{2}\left( u_{1}-u_{2}\right) \right\Vert _{X_{p}}, 
\]%
where $C$ is the constant in Lemma $3.1$. In a similar way for $1\leq k\leq
s $ we obtain 
\begin{equation}
\left\Vert D_{\xi }^{k}\left[ \hat{F}\left( u_{1}\right) -\hat{F}\left(
u_{2}\right) \right] \right\Vert _{X_{p}}\leq  \tag{3.13}
\end{equation}%
\[
M_{0}\bar{f}\left( M+1\right) \left\Vert u_{1}-u_{2}\right\Vert _{X_{\infty
}}+M_{1}\bar{f}\left( M+1\right) \left\Vert D_{\xi }^{k}\left(
u_{1}-u_{2}\right) \right\Vert _{X_{p}}. 
\]%
From $\left( 3.10\right) -\left( 3.13\right) $, using Minkowski's inequality
for integrals, Fourier multiplier theorems for operator-valued functions in $%
X_{p}$ and Young's inequality, we obtain%
\[
\left\Vert G\left( u_{1}\right) -G\left( u_{2}\right) \right\Vert _{Y\left(
T\right) }\leq \left( 2\pi \right) ^{-\frac{1}{n}}\dint\limits_{R^{n}}\dint%
\limits_{0}^{t}\left\Vert u_{1}-u_{2}\right\Vert _{X_{\infty }}+\left\Vert
u_{1}-u_{2}\right\Vert _{Y^{2,p}}d\xi d\tau + 
\]%
\[
\dint\limits_{0}^{t}\left\Vert F\left( u_{1}\right) -F\left( u_{2}\right)
\right\Vert _{X_{\infty }}d\tau +\left( 2\pi \right) ^{-\frac{1}{n}%
}\dint\limits_{0}^{t}\dint\limits_{R^{n}}\left\Vert F\left( u_{1}\right)
-F\left( u_{2}\right) \right\Vert _{Y^{2,p}}d\xi d\tau \leq 
\]%
\[
T\left[ 1+C_{1}\left( M+1\right) ^{2}\bar{f}\left( M+1\right) \right]
\left\Vert u_{1}-u_{2}\right\Vert _{Y\left( T\right) }, 
\]%
where $C_{1}$ is a constant. If $T$ satisfies $\left( 3.9\right) $ and the
following inequality 
\begin{equation}
T\leq \frac{1}{2}\left[ 1+C_{1}\left( M+1\right) ^{2}\bar{f}\left(
M+1\right) \right] ^{-1},  \tag{3.14}
\end{equation}%
then 
\[
\left\Vert Gu_{1}-Gu_{2}\right\Vert _{Y\left( T\right) }\leq \frac{1}{2}%
\left\Vert u_{1}-u_{2}\right\Vert _{Y\left( T\right) }. 
\]%
That is, $G$ is a constructive map. By contraction mapping principle, we
know that $G(u)$ has a fixed point $u(x,t)\in $ $Q\left( M;T\right) $ that
is a solution of the problem $(1.1)-(1.2)$. From $\left( 3.2\right) $ we get
that $u$ is a solution of $\left( 3.5\right) .$ Let us show that this
solution is a unique in $Y\left( T\right) $. Let $u_{1}$, $u_{2}\in Y\left(
T\right) $ are two solution of the problem $(1.1)-(1.2)$. Then%
\begin{equation}
u_{1}-u_{2}=  \tag{3.15}
\end{equation}

\[
\left( 2\pi \right) ^{-\frac{1}{n}}\dint\limits_{R^{n}}\dint\limits_{0}^{t}%
\breve{S}\left( t-\tau ,\xi ,A\right) \left[ \hat{F}\left( u_{1}\right)
\left( \xi ,\tau \right) -\hat{F}\left( u_{2}\right) \left( \xi ,\tau
\right) \right] d\tau d\xi . 
\]%
By the definition of the space $Y\left( T\right) $, we can assume that%
\[
\left\Vert u_{1}\right\Vert _{X_{\infty }}\leq C_{1}\left( T\right) ,\text{ }%
\left\Vert u_{1}\right\Vert _{X_{\infty }}\leq C_{1}\left( T\right) . 
\]%
Hence, by Lemmas 3.1, 3.2, Theorem 2.1 and Minkowski's inequality for
integrals we obtain from $\left( 3.15\right) $

\begin{equation}
\left\Vert u_{1}-u_{2}\right\Vert _{Y^{s,p}}\leq C_{2}\left( T\right) \text{ 
}\dint\limits_{0}^{t}\left\Vert u_{1}-u_{2}\right\Vert _{Y^{s,p}}d\tau . 
\tag{3.16}
\end{equation}%
From $(3.16)$ and Gronwall's inequality, we have $\left\Vert
u_{1}-u_{2}\right\Vert _{Y^{2,p}}=0$, i.e. problem $(1.1)-(1.2)$ has a
unique solution which belongs to $Y\left( T\right) .$ That is, we obtain the
first part of the assertion. Now, let $\left[ 0\right. ,\left. T_{0}\right) $
be the maximal time interval of existence for $u\in Y\left( T_{0}\right) $.
It remains only to show that if $(3.0)$ is satisfied, then $T_{0}=\infty $.
Assume contrary that, $\left( 3.0\right) $ holds and $T_{0}<\infty .$ For $%
T^{\prime }\in \left[ 0\right. ,\left. T_{0}\right) ,$ we consider the
following integral equation

\begin{equation}
\upsilon \left( x,t\right) =S_{1}\left( t,A\right) u\left( x,T\right)
+S_{2}\left( t,A\right) u_{t}\left( x,T\right) +\tilde{\Phi}\left(
t,x\right) .  \tag{3.17}
\end{equation}%
By virtue of $(3.5)$, for $T^{\prime }>T$ we have 
\[
\sup_{T^{\prime }\in \left[ 0\right. ,T_{0}\left. {}\right) }\left\Vert
Au\left( .,T^{\prime }\right) \right\Vert _{Y^{s,p}}+\left\Vert Au\left(
.,T^{\prime }\right) \right\Vert _{X_{\infty }}+\left\Vert Au_{t}\left(
.,T^{\prime }\right) \right\Vert _{Y^{s,p}}+ 
\]%
\[
\left\Vert Au_{t}\left( .,T^{\prime }\right) \right\Vert _{X_{\infty
}}<\infty , 
\]%
where%
\[
u\left( x,T^{\prime }\right) =\varphi \left( x\right)
+\dint\limits_{T^{\prime }}^{T}\alpha \left( \sigma \right) u\left( x,\sigma
\right) d\sigma ,\text{ } 
\]

\[
u_{t}\left( x,T^{\prime }\right) =\psi \left( x\right)
+\dint\limits_{T^{\prime }}^{T}\beta \left( \sigma \right) u_{t}\left(
x,\sigma \right) d\sigma . 
\]%
By reasoning as a first part of theorem and by contraction mapping
principle, there is a $T^{\ast }\in \left( 0,T_{0}\right) $ such that for
each $T^{\prime }\in \left[ 0\right. ,\left. T_{0}\right) ,$ the equation $%
\left( 3.17\right) $ has a unique solution $\upsilon \in Y\left( T^{\ast
}\right) .$ The estimates $\left( 3.9\right) $ and $\left( 3.14\right) $
imply that $T^{\ast }$ can be selected independently of $T^{\prime }\in %
\left[ 0\right. ,\left. T_{0}\right) .$ Set $T^{\prime }=T_{0}-\frac{T^{\ast
}}{2}$ and define 
\[
\tilde{u}\left( x,t\right) =\left\{ 
\begin{array}{c}
u\left( x,t\right) ,\text{ }t\in \left[ 0,T^{\prime }\right] \\ 
\upsilon \left( x,t-T^{\prime }\right) \text{, }t\in \left[ T^{\prime
},T_{0}+\frac{T^{\ast }}{2}\right]%
\end{array}%
\right. . 
\]%
By construction $\tilde{u}\left( x,t\right) $ is a solution of the problem $%
(1.1)-(1.2)$ on $\left[ T^{\prime },T_{0}+\frac{T^{\ast }}{2}\right] $ and
in view of local uniqueness, $\tilde{u}\left( x,t\right) $ extends $u.$ This
is against to the maximality of $\left[ 0\right. ,\left. T_{0}\right) $, i.e
we obtain $T_{0}=\infty .$\qquad

\begin{center}
\textbf{4. \ The nonlocal Cauchy problem for the finite and infinite many
system of wave equation }
\end{center}

Consider first, the linear problem $\left( 1.11\right) .$ Let (see $\left[ 
\text{40, \S\ 1.18}\right] $)%
\[
\text{ }l_{q}\left( N\right) =\left\{ \text{ }u=\left\{ u_{j}\right\} ,\text{
}j=1,2,...N,\left\Vert u\right\Vert _{l_{q}\left( N\right) }=\left(
\sum\limits_{j=1}^{N}\left\vert u_{j}\right\vert ^{q}\right) ^{\frac{1}{q}%
}<\infty \right\} , 
\]%
and 
\[
\text{ }l_{q}^{\sigma }\left( N\right) =\left\{ \text{ }\left\Vert
u\right\Vert _{l_{q}^{\sigma }\left( N\right) }=\left(
\sum\limits_{j=1}^{N}2^{\sigma j}u_{j}^{q}\right) ^{\frac{1}{q}}<\infty
\right\} \text{, }\sigma >0. 
\]

Here, 
\[
X_{pq}=L^{p}\left( R^{n};l_{q}\right) ,Y^{s,p,q}=H^{s,p}\left(
R^{n};l_{q}\right) ,Y_{1,0}^{s,p,q,\sigma }=H^{s,p}\left(
R^{n};l_{q}^{\sigma }\right) \cap L^{1}\left( R^{n};l_{q}^{\sigma }\right) , 
\]

\[
Y_{1,1}^{s,p,q,\sigma }=H^{s,p}\left( R^{n};l_{q}^{\frac{\sigma }{2}}\right)
\cap L^{1}\left( R^{n};l_{q}^{\frac{\sigma }{2}}\right) \text{, }Y_{\infty
}^{s,p,q,\sigma }=H^{s,p}\left( R^{n};l_{q}\right) \cap L^{\infty }\left(
R^{n};l_{q}\right) , 
\]%
\[
\text{ }Y^{s,p,q,\sigma }=H^{s,p}\left( R^{n};l_{q}^{\sigma },l_{q}\right) . 
\]

\textbf{Condition 4.1. }Assume:

(1) $\left\vert 1+\dint\limits_{0}^{T}\alpha \left( \sigma \right) \beta
\left( \sigma \right) d\sigma \right\vert >\dint\limits_{0}^{T}\left(
\left\vert \alpha \left( \sigma \right) \right\vert +\left\vert \beta \left(
\sigma \right) \right\vert \right) d\sigma ;$

(2) $s>n\left( \frac{1}{r}+\frac{1}{p}\right) $ for $r\in \left[ 1,2\right]
, $ $p\in \left[ 1,\infty \right] .$

From Theorem 2.1 we obtain

\textbf{Theorem 4.1. }Assume the Condition 4.1 hold and $0\leq \gamma <\frac{%
1}{2}$. Moreover, $\varphi \in Y_{1,0}^{s,p,q,\sigma },$ $\psi \in
Y_{1,1}^{s,p,q,\sigma }$ and $g\left( .,t\right) \in Y_{1}^{s,p,q}$ for $%
t\in \left[ 0,T\right] .$ Then problem $\left( 1.11\right) $\ has a unique
solution $u\left( x,t\right) \in C^{2}\left( \left[ 0,T\right]
;Y^{s,p,q,\sigma }\right) $ and the following estimate holds 
\[
\left\Vert A^{\gamma }u\left( .,t\right) \right\Vert _{X_{\infty
,q}}+\left\Vert A^{\gamma }u_{t}\left( .,t\right) \right\Vert _{X_{\infty
,q}}\leq C\left\{ \left\Vert A\varphi \right\Vert _{Y_{1,0}^{s,p,q,\sigma
}}+\left\Vert A\varphi \right\Vert _{X_{1,q}}\right. + 
\]

\[
\left\Vert A\psi \right\Vert _{Y_{1,1}^{s,p,q,\sigma }}+\left\Vert A\psi
\right\Vert _{X_{1,q}}+\left. \dint\limits_{0}^{t}\left( \left\Vert \Delta
g\left( .,\tau \right) \right\Vert _{Y^{s,p,q}}+\left\Vert \Delta g\left(
.,\tau \right) \right\Vert _{X_{1,q}}\right) d\tau \right\} 
\]%
uniformly in $t\in \left[ 0,T\right] .$

Consider now, the integral Cauchy problem for the following nonlinear system%
\begin{equation}
\left( u_{m}\right) _{tt}-\Delta \left( u_{m}\right)
_{tt}+\sum\limits_{j=1}^{N}a_{mj}u_{j}\left( x,t\right) =F_{m}\left(
u\right) ,\text{ }x\in R^{n},\text{ }t\in \left( 0,T\right) ,  \tag{4.1}
\end{equation}%
\begin{equation}
u_{m}\left( x,0\right) =\varphi \left( x\right) +\dint\limits_{0}^{T}\alpha
\left( \sigma \right) u_{m}\left( x,\sigma \right) d\sigma ,\text{ } 
\tag{4.2}
\end{equation}

\begin{equation}
\left( u_{m}\right) _{t}\left( x,0\right) =\psi \left( x\right)
+\dint\limits_{0}^{T}\beta \left( \sigma \right) u_{t}\left( x,\sigma
\right) d\sigma ,  \tag{4.3}
\end{equation}%
\[
\text{ }m=1,2,...,N,\text{ }N\in \mathbb{N}, 
\]%
where $u=\left( u_{1},u_{2},...,u_{N}\right) ,$ $a_{mj}$ are complex
numbers, $\varphi _{m}\left( x\right) $ and $\psi _{m}\left( x\right) $ are
data functions$.$ Here,%
\[
l_{q}=\text{ }l_{q}\left( N\right) =\left\{ \text{ }u=\left\{ u_{j}\right\} ,%
\text{ }j=1,2,...N,\left\Vert u\right\Vert _{l_{q}\left( N\right) }=\left(
\sum\limits_{j=1}^{N}\left\vert u_{j}\right\vert ^{q}\right) ^{\frac{1}{q}%
}<\infty \right\} , 
\]%
(see $\left[ \text{44, \S\ 1.18}\right] $)$.$ Let $A$ be the operator in $%
l_{q}\left( N\right) $ defined by%
\[
\text{ }A=\left[ a_{mj}\right] \text{, }a_{mj}=g_{m}2^{sj},\text{ }%
m,j=1,2,...,N,\text{ }D\left( A\right) =\text{ }l_{q}^{s}\left( N\right) = 
\]

\[
\left\{ \text{ }u=\left\{ u_{j}\right\} ,\text{ }j=1,2,...N,\left\Vert
u\right\Vert _{l_{q}^{s}\left( N\right) }=\left(
\sum\limits_{j=1}^{N}2^{sj}u_{j}^{q}\right) ^{\frac{1}{q}}<\infty \right\} . 
\]

From Theorem 3.1 we obtain the following result

\textbf{Theorem 4.1. }Assume:

(1) the Condition 4.1 hold and $0\leq \gamma <\frac{1}{2}$.

(2) the function $u\rightarrow $ $F\left( u\right) $: $R^{n}\times \left[ 0,T%
\right] \times E_{0q}\rightarrow l_{q}$ is a measurable function in $\left(
x,t\right) \in R^{n}\times \left[ 0,T\right] $ for $u\in E_{0q};$

(3) the function $u\rightarrow $ $F\left( u\right) $ is continuous in $u\in
E_{0q}$ for $x,$ $t\in R^{n}\times \left[ 0,T\right] ;$ moreover $F\left(
u\right) \in C^{\left( 1\right) }\left( E_{0q};l_{q}\right) $.

Then problem $\left( 4.1\right) -\left( 4.2\right) $ has a unique local
strange solution%
\[
u\in C^{\left( 2\right) }\left( \left[ 0\right. ,\left. T_{0}\right)
;Y_{\infty }^{s,p,q,\sigma }\right) , 
\]%
where $T_{0}$ is a maximal time interval that is appropriately small
relative to $M$. Moreover, if

\begin{equation}
\sup_{t\in \left[ 0\right. ,\left. T_{0}\right) }\left( \left\Vert
Au\right\Vert _{Y^{s,p,q}}+\left\Vert Au\right\Vert _{X_{\infty
,q}}+\left\Vert Au_{t}\right\Vert _{Y^{s,p,q}}+\left\Vert Au_{t}\right\Vert
_{X_{\infty ,q}}\right) <\infty  \tag{4.4}
\end{equation}%
then $T_{0}=\infty .$

\ \textbf{Proof. }By virtue of $\left[ 43\right] ,$ $l_{q}\left( N\right) $
is a Fourier type space. It is easy to see that the operator $A$ is positive
in $l_{q}\left( N\right) .$ Moreover, by interpolation theory of Banach
spaces $\left[ \text{40, \S\ 1.3}\right] $, we have 
\[
E_{0q}=\left( W^{2,p}\left( R^{n};l_{q}^{s},l_{q}\right) ,L^{p}\left(
R^{n};l_{q}\right) \right) _{\frac{1}{2p},q}=B_{p,q}^{2\left( 1-\frac{1}{2p}%
\right) }\left( R^{n};l_{q}^{s\left( 1-\frac{1}{2p}\right) },l_{q}\right) . 
\]%
By using the properties of spaces $Y^{s,p,q},$ $Y_{\infty }^{s,p,q},$ $%
E_{0q} $ we get that all conditions of Theorem 3.1 are hold, i,e., we obtain
the conclusion.\ 

\begin{center}
\textbf{5.} \textbf{The Wentzell-Robin type mixed problem for wave equations}
\end{center}

Consider at first, the linear problem $\left( 1.7\right) -\left( 1.9\right) $%
. Here, 
\[
X_{p2}=L^{p}\left( R^{n};L^{2}\left( 0,1\right) \right)
,Y^{s,p,2}=H^{s,p}\left( R^{n};L^{2}\left( 0,1\right) \right)
,Y_{1}^{s,p,2}= 
\]%
\[
H^{s,p}\left( R^{n};L^{2}\left( 0,1\right) \right) \cap L^{1}\left(
R^{n};L^{2}\left( 0,1\right) \right) ,\text{ }Y_{\infty
}^{s,p,2}=H^{s,p}\left( R^{n};L^{2}\left( 0,1\right) \right) \cap 
\]

\[
L^{\infty }\left( R^{n};L^{2}\left( 0,1\right) \right) ,\text{ }%
Y^{s,p,2,2}=H^{s,p}\left( R^{n};H^{2}\left( 0,1\right) ,L^{2}\left(
0,1\right) \right) . 
\]

From Theorem 2.1 we obtain the following result

\textbf{Theorem 5.1.} Suppose the the following conditions are satisfied:

(1) The Condition 4.1 is hold and $0\leq \gamma <\frac{1}{2}$;

(2)\ $a$ is positive, $b$ is a real-valued functions on $\left( 0,1\right) $%
. Moreover$,$ $a\left( .\right) \in C\left( 0,1\right) $ and%
\[
\exp \left( -\dint\limits_{\frac{1}{2}}^{x}b\left( t\right) a^{-1}\left(
t\right) dt\right) \in L_{1}\left( 0,1\right) ; 
\]

(3) $\varphi ,$ $\psi \in Y_{1}^{s,p}$ and $g\left( .,t\right) \in
Y_{1}^{s,p}$ for $t\in \left[ 0,T\right] .$

Then the problem $\left( 1.7\right) -\left( 1.9\right) $\ has a unique
solution $u\left( x,t\right) \in C^{2}\left( \left[ 0,T\right]
;Y^{s,p,2}\right) $ and the following estimate holds 
\begin{equation}
\left\Vert A^{\gamma }u\right\Vert _{X_{\infty },2}+\left\Vert A^{\gamma
}u_{t}\right\Vert _{X_{\infty },2}\leq C\left\{ \left\Vert A\varphi
\right\Vert _{Y^{s,p,2}}+\left\Vert A\varphi \right\Vert _{X_{1},2}\right. +
\tag{5.1}
\end{equation}

\[
\left\Vert A\psi \right\Vert _{Y^{s,p},2}+\left\Vert A\psi \right\Vert
_{X_{1},2}+\left. \dint\limits_{0}^{t}\left( \left\Vert \Delta g\left(
.,\tau \right) \right\Vert _{Y^{s,p,2}}+\left\Vert \Delta g\left( .,\tau
\right) \right\Vert _{X_{1},2}\right) d\tau \right\} 
\]%
uniformly in $t\in \left[ 0,T\right] .$

\ \textbf{Proof.} Let $H=L^{2}\left( 0,1\right) $ and $A$ is a operator
defined by $\left( 1.6\right) .$ Then the problem $\left( 1.7\right) -\left(
1.9\right) $ can be rewritten as the problem $\left( 2.1\right) -\left(
2.2\right) $. By virtue of $\left[ \text{12, 25}\right] $ the operator $A$
generates analytic semigroup in $L^{2}\left( 0,1\right) $. Hence, by virtue
of (1)-(3) all conditions of Theorem 2.1 are satisfied. Then Theorem 2.1
implies the assertion.

Consider now, the nonlinear problem 
\begin{equation}
u_{tt}-\Delta u+au_{yy}+bu_{y}+F\left( u\right) =0,\text{ }  \tag{5.2}
\end{equation}%
\ \ \ 

\begin{equation}
u\left( x,y,0\right) =\varphi \left( x,y\right) +\dint\limits_{0}^{T}\alpha
\left( \sigma \right) u\left( x,y,\sigma \right) d\sigma ,\text{ }  \tag{5.3}
\end{equation}

\[
u_{t}\left( x,y,0\right) =\psi \left( x,y\right) +\dint\limits_{0}^{T}\beta
\left( \sigma \right) u_{t}\left( x,y,\sigma \right) d\sigma , 
\]

\begin{equation}
a\left( j\right) u_{yy}\left( x,j,t\right) +b\left( j\right) u_{y}\left(
x,j,t\right) =0,\text{ }j=0,1\text{ for all }t\in \left[ 0,T\right] . 
\tag{5.4}
\end{equation}

\bigskip Let%
\[
\text{ }E_{0,2}=\left( X_{p2},Y^{s,p,2,2}\right) _{\frac{1}{2p},2}. 
\]

Now, we obtain the following result:

\bigskip \textbf{Theorem 5.2. } Suppose the the following conditions are
satisfied:

1.\ $a$ is positive, $b$ is a real-valued functions on $\left( 0,1\right) $.
Moreover$,$ $a\left( .\right) \in C\left( 0,1\right) $ and%
\[
\exp \left( -\dint\limits_{\frac{1}{2}}^{x}b\left( t\right) a^{-1}\left(
t\right) dt\right) \in L_{1}\left( 0,1\right) ; 
\]

2. the Condition 4.1 is hold and $0\leq \gamma <1$;

3. (4) and (5) assumptions of the Conditions 3.1 are satisfied for $%
E=L^{2}\left( 0,1\right) $;

4. the function%
\[
u\rightarrow F\left( u\right) :R^{n}\times \left[ 0,T\right] \times
E_{02}\rightarrow L^{2}\left( 0,1\right) 
\]%
is a measurable in $\left( x,t\right) \in R^{n}\times \left[ 0,T\right] $
for $u\in E_{02};$

5. the function $u\rightarrow $ $F\left( u\right) $ is continuous in $u\in
E_{02}$ for $x,$ $t\in R^{n}\times \left[ 0,T\right] ;$ moreover $f\left(
u\right) \in C^{\left( 1\right) }\left( E_{02};L^{2}\left( 0,1\right)
\right) $.

Then problem $\left( 5.2\right) -\left( 5.4\right) $ has a unique local
strange solution%
\[
u\in C^{\left( 2\right) }\left( \left[ 0\right. ,\left. T_{0}\right)
;Y_{\infty }^{s,p,2,2}\right) , 
\]%
where $T_{0}$ is a maximal time interval that is appropriately small
relative to $M$. Moreover, if

\begin{equation}
\sup_{t\in \left[ 0\right. ,\left. T_{0}\right) }\left( \left\Vert
Au\right\Vert _{Y^{s,p,2}}+\left\Vert Au\right\Vert _{X_{\infty
,2}}+\left\Vert Au_{t}\right\Vert _{Y^{s,p,2}}+\left\Vert Au_{t}\right\Vert
_{X_{\infty ,2}}\right) <\infty  \tag{5.5}
\end{equation}%
then $T_{0}=\infty .$

\ \textbf{Proof. }By virtue of $\left[ 43\right] ,$ $L^{2}\left( 0,1\right) $
is a Fourier type space. It is easy to see that the operator $A$ is positive
in $L^{2}\left( 0,1\right) .$ Moreover, by interpolation theory of Banach
spaces $\left[ \text{40, \S\ 1.3}\right] $, we have 
\[
E_{02}=\left( W^{s,p}\left( R^{n};H^{2}\left( 0,1\right) ,L^{2}\left(
0,1\right) \right) ,L^{p}\left( R^{n};L^{2}\left( 0,1\right) \right) \right)
_{\frac{1}{2p},2}= 
\]%
\[
B_{p,2}^{s\left( 1-\frac{1}{2p}\right) }\left( R^{n};H^{s\left( 1-\frac{1}{2p%
}\right) },L^{2}\left( 0,1\right) \right) . 
\]%
By using the properties of spaces $Y^{s,p,2},$ $Y_{\infty }^{s,p,2},$ $%
E_{02} $ we get that all conditions of Theorem 3.1 are hold, i,e., we obtain
the conclusion.\ 

\bigskip \textbf{References}

\begin{quote}
\ \ \ \ \ \ \ \ \ \ \ \ \ \ \ \ \ \ \ \ \ \ \ \ 
\end{quote}

\begin{enumerate}
\item H. Amann, Linear and quasi-linear equations,1, Birkhauser, Basel
1995.\ \ 

\item A. Arosio, Linear second order differential equations in Hilbert
space. The Cauchy problem and asymptotic behaviour for large time, Arch.
Rational Mech. Anal., v 86, (2), 1984, 147-180.

\item A. Ashyralyev, N. Aggez, Nonlocal boundary value hyperbolic problems
Involving Integral conditions, Bound.Value Probl., 2014 V. 2014:214.

\item J. Bourgain, A Hausdor -Young inequality for B-convex Banach spaces,
Paci c J. Math. 10(1982), 2, 255-262.

\item H. Bahouri, P. G\'{e}rard, High frequency approximation of solutions
to critical nonlinear wave equations, Amer. J. Math. 121 (1999), 131--175.

\item J. Colliander, M. Keel, G. Staffilani, H. Takaoka, T. Tao, Sharp
global well-posedness for KdV and modified KdV on R and T, J. Amer. Math.
Soc. 16 (2003), 705--749.

\item G. Da Prato and E. Giusti, A characterization on generators of
abstract cosine functions,\ Boll. Del. Unione Mat., (22)1967, 367--362.

\item Denk R., Hieber M., Pr\"{u}ss J., $R$-boundedness, Fourier multipliers
and problems of elliptic and parabolic type, Mem. Amer. Math. Soc. 166
(2003), n.788.

\item E. B. Davies and M. M. Pang, The Cauchy problem and a generalization
of the Hille-Yosida theorem, Proc. London Math. Soc. (55) 1987, 181-208.

\item D. Foschi, S. Klainerman, Bilinear space-time estimates for
homogeneous wave equations, Ann. Sci. de l'ENS, 33 no.2 (2000).

\item H. O. Fattorini, Second order linear differential equations in Banach
spaces, in North Holland Mathematics Studies,\ V. 108, North-Holland,
Amsterdam, 1985.

\item A. Favini, G. R. Goldstein, J. A. Goldstein and S. Romanelli,
Degenerate second order differential operators generating analytic
semigroups in $L_{p}$ and $W^{1,p}$, Math. Nachr. 238 (2002), 78 --102.

\item A. Favini, V. Shakhmurov, Y. Yakubov, Regular boundary value problems
for complete second order elliptic differential-operator equations in UMD
Banach spaces, Semigroup form, v. 79 (1), 2009.

\item J. Ginibre, G. Velo, Generalized Strichartz Inequalities for the Wave
Equation, Jour. Func. Anal., 133 (1995), 50--68.

\item M. Girardy, L. Weis, Operator-valued multiplier theorems on Besov
spaces, Math. Nachr. 251 (2003), 34-51.

\item J. A. Goldstein, Semigroup of linear operators and applications,
Oxford, 1985.

\item V. I. Gorbachuk and M. L. Gorbachuk M, Boundary value problems \ for
differential-operator equations, Naukova Dumka, Kiev, 1984.

\item M. Grillakis, Regularity for the wave equation with a critical
nonlinearity, Comm. Pure Appl. Math., 45 (1992), 749--774.

\item S. G. Krein, Linear differential equations in Banach space,
Providence, 1971.

\item C. Kenig, G. Ponce, L. Vega, A bilinear estimate with applications to
the KdV equation, J. Amer. Math. Soc. 9 (1996), 573--603.

\item L. Kapitanski, Weak and Yet Weaker Solutions of Semilinear Wave
Equations, Comm. Part. Diff. Eq., 19 (1994), 1629--1676.

\item S. Klainerman, Global existence for nonlinear wave equations, Comm.
Pure Appl. Math. 33 (1980), 43--101.

\item S. Klainerman, M. Machedon, On the regularity properties of a model
problem related to wave maps, Duke Math. J. 87 (1997), 553--589.

\item M. Keel, T. Tao, Local and global well-posedness of wave maps on R1+1
for rough data, IMRN 21 (1998), 1117--1156.

\item V. Keyantuo, M. Warma, The wave equation with Wentzell--Robin boundary
conditions on Lp-spaces, J. Differential Equations 229 (2006) 680--697.

\item H. Lindblad, C. D. Sogge, Restriction theorems and semilinear
Klein-Gordon equations in (1+3) dimensions, Duke Math. J. 85 (1996), no. 1,
227--252.

\item J. L. Lions and J. Peetre, Sur une classe d'espaces d'interpolation,
Inst. Hautes Etudes Sci. Publ. Math., 19(1964), 5-68.

\item A. Lunardi, Analytic semigroups and optimal regularity in parabolic
problems, Birkhauser, 2003.

\item V. G. Makhankov, Dynamics of classical solutions (in non-integrable
systems), Phys. Lett. C 35, (1978) 1--128.

\item N. Masmoudi, K. Nakanishi, From nonlinear Klein-Gordon equation to a
system of coupled nonlinear Schrdinger equations, Math. Ann. 324 (2002), (2)
359--389.

\item L. Nirenberg, On elliptic partial differential equations, Ann. Scuola
Norm. Sup. Pisa 13 (1959), 115--162.

\item S. Piskarev and S.-Y. Shaw, Multiplicative perturbations of semigroups
and applications to step responses and cumulative outputs, J. Funct. Anal.
128 (1995), 315-340.

\item A. Pazy, Semigroups of linear operators and applications to partial
differential equations. Springer, Berlin, 1983.

\item C. D. Sogge, Lectures on Nonlinear Wave Equations, International
Press, Cambridge, MA, 1995.

\item V. B. Shakhmurov, Embedding and separable differential operators in
Sobolev-Lions type spaces, Mathematical Notes, 2008, v. 84, no 6, 906-926.

\item V. B. Shakhmurov, Nonlinear abstract boundary value problems in
vector-valued function spaces and applications, Nonlinear Anal-Theor., v.
67(3) 2006, 745-762.

\item V. B. Shakhmurov, Embedding theorems in Banach-valued B-spaces and
maximal B-regular differential operator equations, Journal of Inequalities
and Applications, v. 2006, 1-22.

\item R. Shahmurov, On strong solutions of a Robin problem modeling heat
conduction in materials with corroded boundary, Nonlinear Anal., Real World
Appl., v.13, (1), 2011, 441-451.

\item R. Shahmurov, Solution of the Dirichlet and Neumann problems for a
modified Helmholtz equation in Besov spaces on an annuals, J. Differential
equations, v. 249(3), 2010, 526-550.

\item H. Triebel, Interpolation theory, Function spaces, Differential
operators, North-Holland, Amsterdam, 1978.

\item H. Triebel, Fractals and spectra, Birkhauser Verlag, Related to
Fourier analysis and function spaces, Basel, 1997.

\item T. Tao, \ Local and global analysis of nonlinear dispersive and wave
equations, CBMS regional conference series in mathematics, 2006.

\item T. Tao, Low regularity semi-linear wave equations. Comm. Partial
Differential Equations, 24 (1999), no. 3-4, 599--629.

\item C. Travis and G. F. Webb, Second order differential equations in
Banach spaces, Nonlinear Equations in Abstract Spaces, (ed. by V.
Lakshmikantham), Academic Press, 1978, 331-361.

\item G. B. Whitham, Linear and Nonlinear Waves, Wiley--Interscience, New
York, 1975.

\item S. Wang, G. Chen,The Cauchy Problem for the Generalized IMBq equation
in $W^{s,p}\left( R^{n}\right) $, J. Math. Anal. Appl. 266, 38--54 (2002).

\item S. Yakubov and Ya. Yakubov, Differential-operator Equations. Ordinary
and Partial \ Differential Equations , Chapman and Hall /CRC, Boca Raton,
2000.

\item N.J. Zabusky, Nonlinear Partial Differential Equations, Academic
Press, New York, 1967.
\end{enumerate}

\end{document}